\documentclass[11pt]{amsart}
\usepackage{etex}

\usepackage{amsmath}
\usepackage{cases}
\usepackage{mathrsfs, euscript}
\usepackage{bbm}
\usepackage{amssymb}
\usepackage{txfonts}
\usepackage{amscd}
\usepackage{amsfonts,latexsym,amsmath,amsthm,amsxtra,mathdots,amssymb,latexsym}
\usepackage[all,cmtip]{xy}
\usepackage{pictexwd, amsthm, amssymb, amsmath, amsxtra, graphicx, psfrag, latexsym}

\usepackage{enumerate}
\RequirePackage{amsmath} \RequirePackage{amssymb}
\usepackage{color}
\usepackage{colordvi}
\usepackage{multicol}
\usepackage{hyperref}
\usepackage{footnote}

    \newcommand{\Isom}{{\mathrm{Isom}}}

     \newcommand{\rank}{{\mathrm{rank}}}

\def\-{^{-1}}

\def\-{^{-1}}

\newcommand{\delete}[1]{}

     \newcommand{\SL}{{\mathrm{SL}}}
     \newcommand{\SO}{{\mathrm{SO}}}
    \newcommand{\SU}{{\mathrm{SU}}}

        \newcommand{\Sp}{{\mathrm{Sp}}}

    \theoremstyle{plain}

\newtheorem*{main}{Main Theorem}

    \newtheorem{thm}{Theorem}[section] \newtheorem{cor}[thm]{Corollary}
    \newtheorem{lem}[thm]{Lemma}  \newtheorem{prop}[thm]{Proposition}
    
    \newtheorem {conj}[thm]{Conjecture} \newtheorem{defn}[thm]{Definition}

    \numberwithin{equation}{section}

\theoremstyle{remark}
\newtheorem*{remark}{Remark}

\begin{document}

\title{On splitting rank of non-compact type symmetric spaces and bounded cohomology}

\author[Shi Wang]{Shi Wang}
\address{Department of Mathematics,
Indiana University,
831 E. Third St.,
Bloomington, IN 47405, U.S.A.}
\email{wang679@iu.edu}

\begin{abstract}
Let $X=G/K$ be a higher rank symmetric space of non-compact type where $G=\Isom^0(X)$. We define the splitting rank of $X$, denoted by $\text{srk}(X)$, to be the maximal dimension of a totally geodesic submanifold $Y\subset X$ which splits off an isometric $\mathbb R$-factor. We compute explicitly the splitting rank for each irreducible symmetric space. For an arbitrary (not necessarily irreducible) symmetric space, we show that the comparison map $\eta:H^{*}_{c,b}(G,\mathbb{R})\rightarrow H^{*}_c(G,\mathbb{R})$ is surjective in degrees $*\geq \text{srk}(X)+2$, provided $X$ has no direct factors of $\mathbb H^2$, $\SL(3,\mathbb R)/\SO(3)$, $\Sp(2,\mathbb R)/U(2)$, $G_2^2/\SO(4)$ and $\SL(4,\mathbb R)/\SO(4)$. This generalizes the result of \cite{LW} regarding Dupont's problem.
\end{abstract}

\maketitle


\section{Introduction}

The notion of bounded cohomology was introduced by Gromov (see \cite{Gr}). Let $X$ be a topological space, we denote $C^n(X)$ the group of singular $n$-cochains with real coefficients. We say an element $c$ in this cochain group is bounded, if the values of $c$ on each $n$-simplex is bounded, that is, the norm $||c||=sup\,\{\,|c(\sigma)|:\sigma:\Delta^n\rightarrow X \;\textrm{is continuous}\,\}$ is finite. We denote $C_b^n(X)$ the set of all bounded $n$-cochains, and the corresponding chain complex induces the bounded singular cohomology $H^n_b(X)$. The inclusion map $i:C_b^n(X)\rightarrow C^n(X)$ induces the comparison map $\phi:H^n_b(X) \rightarrow H^n(X)$.

This notion of bounded cohomology can be extended to groups. For a discrete group $\Gamma$, the group cohomology with real coefficients $H^n(\Gamma)$ can be defined by the cochain complex $C^n(\Gamma)=\{f:\Gamma^n\rightarrow \mathbb{R}\}$ together with a certain coboundary operator. Here we are using the inhomogeneous complex, but one can also work on a homogeneous complex that takes all group invariant functions on $(n+1)$-tuples, yet with a different coboundary operator. Similarly, one defines the bounded group cohomology $H^n_b(\Gamma)$ using a subcomplex, namely the bounded cochains $C^n_b(\Gamma)=\{f:\Gamma^n\rightarrow \mathbb{R}\mid f\,\textrm{is bounded}\}$. For a topological group $G$, so as not to lose its topological information, we consider the continuous cochain complex $C^n_c(G)=\{f:G^n\rightarrow \mathbb{R}\mid f\,\textrm{is continuous}\}$ and correspondingly the continuous bounded cochains $C^n_{c,b}(G)=\{f:G^n\rightarrow \mathbb{R}\mid f\,\textrm{is continuous and bounded}\}$. This gives rise to the continuous cohomology $H^n_c(G)$ and the continuous bounded cohomology $H^n_{c,b}(G)$.

It is worth pointing out that all these notions can be defined in a most generalized way (with arbitrary coefficients and with arbitrary group representations). However in this article we only focus on trivial real coefficients, as is defined above, and meanwhile we refer the readers to \cite{Jo}, \cite{Mo2}, \cite{Ro} for their interests.

Despite the fact that bounded cohomology is easily defined, little is known about these groups in general. It was shown by Gromov \cite{Gr} that $H^*_b(\pi_1(M))\simeq H^*_b(M)$ through the classifying map, so one might just focus on the study of the bounded cohomology for groups. It is clear that $H^0_b(\Gamma)\simeq \mathbb{R}$ and $H^1_b(\Gamma)$ vanishes following the definition, and it was pointed out in \cite{Gr} that $H^n_b(\Gamma)=0\;( n\geq 1)$ for any amenable group $\Gamma$ (following work of Hirsch and Thurston \cite{HT}). However, the bounded cohomology is hard to compute in general--there is only known results for specific groups and in specific (often low) degrees. See \cite{Br}, \cite{Gri}, \cite{So} for free group $\mathbb{F}_2$.

One way to study the bounded cohomology is to look at the comparison map, from the bounded cohomology to the ordinary cohomology, and it is natural to ask whether this map is surjective. The answer is of course no in general. For example, we can easily construct an abelian group $\mathbb{Z}\oplus \mathbb{Z}$ where the ordinary cohomology is nonvanishing in degree two, but the bounded cohomology vanishes due to amenability. However, surjectivity might still hold for certain classes of groups. For example, in the case of semisimple Lie groups, Dupont \cite{Du2} (see also \cite[Problem A']{Mo}, and \cite[Conjecture 18.1]{BIMW}) conjectured the following:

\begin{conj} Let $G$ be a connected semisimple Lie group with finite center, the comparison map
$$H^{*}_{c,b}(G,\mathbb{R})\rightarrow H^{*}_c(G,\mathbb{R})$$
is always surjective.
\end{conj}

This conjecture remains open in the specific case of $\SL(n,\mathbb{R})$. Prior work includes Hartnick and Ott \cite{HO}, which confirmed the conjecture for Lie groups of Hermitian type (as well as some other cases). Domic and Toledo gave explicit bounds in degree two \cite{DT}, and this was later generalized by Clerc and {\O}rsted in \cite{CO}. Lafont and Schmidt \cite{LS} showed surjectivity on top degree (the dimension of the corresponding symmetric space) in all cases excluding $\SL(3,\mathbb{R})$, followed by Bucher-Karlsson's complementary result \cite{Bu}, thus completing an equivalent conjecture of Gromov: that the simplicial volume of any closed locally symmetric space of noncompact type is positive. One of the key step in their approach is to show boundedness of a certain Jacobian, which relies heavily on previous work of Connell and Farb \cite{CF1}, \cite{CF2}. Recently, Inkang Kim and Sungwoon Kim \cite{KK} extended the Jacobian estimate to codimension one (but the codimension one surjectivity of the comparison map is automatic), and they also gave detailed investigation on rank two cases. Meanwhile, Lafont and Wang \cite{LW} showed surjectivity in codimesion $\leq \rank(X)-2$, in irreducible cases excluding $\SL(3,\mathbb R)$ and $\SL(4,\mathbb R)$. In this paper, we extend their results to smaller degrees and show the following:

\begin{main}
Let $X=G/K$ be an $n$-dimensional symmetric space of non-compact type of rank $r\geq 2$, and $\Gamma$ a cocompact torsion-free lattice in $G$. Assume $X$ has no direct factors of $\mathbb H^2$, $\SL(3,\mathbb R)/\SO(3)$, $\Sp(2,\mathbb R)/U(2)$, $G_2^2/\SO(4)$ and $\SL(4,\mathbb R)/\SO(4)$, then the comparison maps
$\eta:H^{*}_{c,b}(G,\mathbb{R})\rightarrow H^{*}_c(G,\mathbb{R})$ and $\eta':H^{*}_{b}(\Gamma,\mathbb{R})\rightarrow
H^{*}(\Gamma,\mathbb{R})$ are both surjective in all degrees $* \geq \text{srk}(X)+2$.
\end{main}

\begin{remark} We will see below in Corollary \ref{cor:si>rank} that $\text{srk}(X)$ is in general smaller than $(n-r)$ (unless the Lie group is of type $\SL(r+1,\mathbb R)$ that makes them equal, in which case we recover the main theorem of \cite{LW}). This means that our Main Theorem generically provides larger range where the comparison map is surjective. On the other hand, the method of barycentric straightening that we are using fails at the splitting rank (See \cite[Theorem 5.6]{LW}, and see Table 1 for explicit expressions), that is, we shall not expect a proof of surjectivity in degrees $\leq\text{srk}(X)$ via this method. In this sense, our main theorem almost fills in the gap (between $\text{srk}(X)$ and $n-r+2$), while leaving the only unknown degree at $\text{srk}(X)+1$.

\end{remark}

{\bf Example of $\SL(4,\mathbb C):$}
If $X=\SL(4,\mathbb C)/\SU(4)$, then $\dim(X)=15$, and according to Table 1 below $\text{srk}(X)=9$. The continuous cohomology of $H_c^*(\SL(4,\mathbb C))$ is isomorphic to the cohomology of its compact dual symmetric space $H^*(\SU(4))$, that is, the exterior algebra generated by $\{\alpha_3, \alpha_5, \alpha_7\}$. According to our main theorem, all classes are bounded if the degree is at least $11$, thus we obtain that the class $\alpha_5\wedge\alpha_7$ is bounded (the top class $\alpha_3\wedge\alpha_5\wedge\alpha_7$ is previously known to be bounded).

\begin{remark}
 Our approach uses the same kind of machinery as \cite{LW}, which will be discussed in details in Section \ref{sec:BC}. Essentially the proof follows similarly except that \cite[Lemma 4.6]{LW} is now replaced by a stronger result (Theorem \ref{thm:brain}). The rest of the paper is devoted to show this theorem hence improving the bounds, but showing this theorem requires an estimate on the $k$-th splitting rank that needs further working case by case on all irreducible symmetric spaces of non-compact type. For the convenience, we leave the most of the computations and analysis in the appendix.
\end{remark}


\section{Splitting rank}\label{SR}
The notion of splitting rank was defined in \cite{LW}, as an obstruction in degree to a certain type of Jacobian being uniformly bounded (See \cite[Theorem 5.6]{LW}). In section \ref{sec:srk}, we will compute in Table 1 the splitting rank of all irreducible symmetric spaces of non-compact type. In section \ref{sec:srk^k}, we will define and analyze the $k$-th splitting rank. And finally in section \ref{sec:reducible}, we will generalize to the reducible cases and establish Theorem \ref{thm:brain}, to give a key estimate on the $k$-th splitting rank.

\subsection{Totally geodesic submanifolds with $\mathbb{R}$-factor}\label{sec:srk} We recall the following definition of splitting rank from \cite{LW}.

\begin{defn}
For $X$ a symmetric space of non-compact type, we define the splitting rank of $X$, denoted $\text{srk}(X)$,
to be the maximal dimension of a totally geodesic submanifold $Y\subset X$ which splits off an isometric
$\mathbb{R}$-factor.
\end{defn}

\begin{remark}
We notice a similar notion of maximal totally geodesic submanifolds is discussed in \cite{BO}, but our definition is slightly different. Indeed, if $X=G_2^2/\SO(4)$, then the submanifold that has dimension equal to the splitting rank is $\mathbb H^2\times \mathbb R$. This is not maximal totally geodesic, since $\mathbb H^2\times \mathbb R\subset \SL(3,\mathbb R)/\SO(3)\subset G_2^2/\SO(4)$ gives a chain of totally geodesic inclusions.
\end{remark}

For a totally geodesic submanifold of a symmetric space, its tangent space can be identified with a Lie triple system. Let $X$ be a symmetric space of non-compact type. We can write $X=G/K$ where $G=\Isom^0(X)$ is the connected component of isometry group of $X$ and $K$ is a maximal compact subgroup of $G$. Fixing a base point $p\in X$, we have the Cartan decomposition $\mathfrak{g}=\mathfrak{k}+\mathfrak{p}$ and $\mathfrak{p}$ can be identified with the tangent space of $X$ at $p$. The following proposition characterizes in terms of Lie algebra, the totally geodesic submanifolds that attain the spitting rank.

\begin{prop}\label{prop:max}
Suppose a totally geodesic submanifold $Y \times \mathbb{R}\subset X$ has dimension equal to the splitting rank of $X$. Then the corresponding Lie triple system $[\mathfrak{p}',[\mathfrak{p}',\mathfrak{p}']]\subset \mathfrak{p}'$ has the form $\mathfrak{p}'=\mathfrak{a}\bigoplus_{\alpha\in \Lambda^+, \alpha(V)=0}\mathfrak{p}_\alpha$, where $\mathfrak{a}$ is a choice of maximal abelian subalgebra in $\mathfrak{p}$ that contains the $\mathbb{R}$-factor $V$.
\end{prop}
\begin{proof}
We identify the tangent space of $X$ with $\mathfrak{p}$ via the Cartan decomposition, and the tangent space of $Y \times \mathbb{R}$ with a Lie triple system $\mathfrak{p}''\subset \mathfrak{p}$. The product structure implies that any vector $v\in \mathfrak{p}''$ commutes with the $\mathbb{R}$-factor $V$. Hence $\mathfrak{p}''\subset \mathfrak{p}'$, where $\mathfrak{p}'=\{\;Z\in \mathfrak{p}\mid [Z,V]=0\;\}$. Notice that $\mathfrak{p}'$ is itself a Lie triple system. To see this, we first extend $V$ to a maximal abelian subalgebra $\mathfrak{a}\subset \mathfrak{p}$, and form the restricted root space decomposition $\mathfrak{p}=\mathfrak{a}\bigoplus_{\alpha\in \Lambda^+} \mathfrak{p}_\alpha$. Then $\mathfrak{p}'$ decomposes as $\mathfrak{a}\bigoplus_{\alpha\in \Lambda^+, \alpha(V)=0}\mathfrak{p}_\alpha$. By a standard Lie algebra computation, we see that $[\mathfrak{p}',\mathfrak{p}']\subset \mathfrak{k}'$, where $\mathfrak{k}'=\mathfrak{k}_0\bigoplus_{\alpha\in \Lambda^+, \alpha(V)=0}\mathfrak{k}_\alpha$, and also $[\mathfrak{k}',\mathfrak{p}']\subset \mathfrak{p}'$. Therefore $\mathfrak{p}'$ is a Lie triple system that contains $\mathfrak{p}''$. By the assumption that $\mathfrak{p}''$ has maximal dimension, we conclude $\mathfrak{p}'=\mathfrak{p}''$. This completes the proof.
\end{proof}

\begin{remark}
We comment that the totally geodesic submanifold in the above proposition is the same as $F(\gamma)$ --the union of all flats that goes through the geodesic $\gamma$ corresponding to the $\mathbb{R}$-factor. In general, $F(\gamma)=F_s(\gamma)\times \mathbb{R}^t$ where $F_s(\gamma)$ is also a symmetric space of non-compact type and $t$ is some integer that measures the singularity of $\gamma$ (see \cite[Proposition 2.20.10]{Eb} for more details). We see in the next proposition that $F(\gamma)$ attains maximal dimension only when $t=1$.
\end{remark}

\begin{prop}
Suppose a totally geodesic submanifold $Y \times \mathbb{R}\subset X$ has dimension equal to the splitting rank of $X$. Then $Y$ is also a symmetric space of non-compact type (i.e. it does not split off an $\mathbb R$-factor).
\end{prop}

\begin{proof}
The proposition is a direct consequence of \cite[Proposition 2.20.10]{Eb} and Proposition \ref{prop:srk gap} below.
\end{proof}

We continue to analyze $Y$ via the above splitting of the Lie algebra. Let $\mathfrak{a}$ be a maximal abelian subalgebra containing the $\mathbb{R}$-factor $V$, and denote by $\mathfrak{a}'\subset \mathfrak{a}$ the orthogonal complement of $V$. Then the Lie triple system of $Y\times \mathbb{R}$ can be written as $V\oplus \mathfrak{a}'\bigoplus_{\alpha\in \Lambda^+, \alpha(V)=0}\mathfrak{p}_\alpha$, where $V$ represents the $\mathbb R$-factor, and $\mathfrak{a}'\bigoplus_{\alpha\in \Lambda^+, \alpha(V)=0}\mathfrak{p}_\alpha$ is the Lie triple system of $Y$. As $Y$ corresponds to a maximal parabolic subalgebra in $\mathfrak{g}$, we can choose a simple system $\Omega=\{\alpha_1,...,\alpha_r\}\subset \Lambda$ corresponding to $X$ such that $\Omega'=\{\alpha_1,...,\alpha_{r-1}\}\subset \ker(V)\cap \Lambda$ is a simple system corresponding to $Y$ (See \cite[Proposition 7.76]{Kn}). In other words, $Y$ has a truncated simple system generated by throwing away one element from the simple system of $X$. We give more detailed information in the next theorem, by simply working through all the cases of irreducible symmetric spaces of non-compact type.

\begin{thm}\label{thm:srk}
Let $X$ be an irreducible symmetric space of non-compact type. Assume $\dim(X)=n$ and $rank(X)=r\geq 2$. We give in the following table the splitting rank of $X$, as well as all totally geodesic submanifolds $Y\times \mathbb{R}$ whose dimension attains the splitting rank.
\end{thm}

\begin{remark}
In the table, we write $\SO^0_{i,j}/\SO_i\SO_j$ short for $\SO_0(i,j)/\SO(i)\times \SO(j)$, and similarly for $\SU_{i,j}/S(U_i U_j)$ and $\Sp_{i,j}/\Sp_i \Sp_j$. We use the same abbreviation in Table 2.
\end{remark}

\renewcommand{\arraystretch}{1.5}
\begin{equation*}
\begin{array}{|c|c|c|c|c|}
\hline
  X & Y & \text{srk}(X) & n & \text{Comments}\\ \hline
  \SL(r+1,\mathbb{R})/\SO(r+1) & \SL(r,\mathbb{R})/\SO(r) & n-r & r(r+3)/2 &r\geq 2\\
  \SL(r+1,\mathbb{C})/\SU(r+1) & \SL(r,\mathbb{C})/\SU(r) & n-2r & r(r+2) &r\geq 2\\
  \SU^*(2r+2)/\Sp(r+1) & \SU^*(2r)/\Sp(r) & n-4r & r(2r+3) &r\geq 2\\
  E_6^{-26}/F_4 & \mathbb{H}^9 & 10 & 26 &r=2\\ \hline
  \SO^0_{r,r+k}/\SO_r\SO_{r+k} &\SO^0_{r-1,r-1+k}/\SO_{r-1} \SO_{r-1+k} & n-(2r+k-2) & r(r+k) &r\geq 2,k\geq 1\\
  \SO(2r+1,\mathbb{C})/\SO(2r+1) & \SO(2r-1,\mathbb{C})/\SO(2r-1) & n-(4r-2) & r(2r+1) &r\geq 2\\\hline
\Sp(r,\mathbb{R})/U(r) & \Sp(r-1,\mathbb{R})/U(r-1) & n-(2r-1) & r(r+1) &r\geq 3 \\
  \SU_{r,r}/S(U_r U_r) & \SU_{r-1,r-1}/S(U_{r-1} U_{r-1}) & n-(4r-3) & 2r^2 &r\geq 3 \\
  \Sp(r,\mathbb{C})/\Sp(r) & \Sp(r-1,\mathbb{C})/\Sp(r-1) & n-(4r-2) & r(2r+1)&r\geq 3  \\
  \SO^*(4r)/U(2r) & \SO^*(4r-4)/U(2r-2) & n-(8r-7) & 2r(2r-1) &r\geq 4\\
  \SO^*(12)/U(6) & \SU^*(6)/\Sp(3) & 15 & 30 &r= 3\\
  \Sp_{r,r}/\Sp_r \Sp_r & \Sp_{r-1,r-1}/\Sp_{r-1} \Sp_{r-1} & n-(8r-5) & 4r^2& r\geq 2\\
  E_7^{-25}/E_6\times U(1) & E_6^{-26}/F_4 & 27 & 54 &r=3\\ \hline
  \SO^0_{r,r}/\SO_r \SO_r & \SO^0_{r-1,r-1}/\SO_{r-1}\SO_{r-1} & n-(2r-2) & r^2 &r\geq 4\\
  \SO(2r,\mathbb{C})/\SO(2r) & \SO(2r-2,\mathbb{C})/\SO(2r-2) & n-(4r-4) & r(2r-1)  &r\geq 4\\ \hline
  \SU_{r,r+k}/S(U_r U_{r+k}) & \SU_{r-1,r-1+k}/S(U_{r-1} U_{r-1+k}) & n-(4r+2k-3) & 2r(r+k) &r\geq 1,k\geq 1\\
  \Sp_{r,r+k}/\Sp_r \Sp_{r+k} & \Sp_{r-1,r-1+k}/\Sp_{r-1} \Sp_{r-1+k} & n-(8r+4k-5) & 4r(r+k) &r\geq 1,k\geq 1\\
  \SO^*(4r+2)/U(2r+1) & \SO^*(4r-2)/U(2r-1) & n-(8r-3) & 2r(2r+1) &r\geq 2\\
  E_6^{-14}/\text{Spin}(10)\times U(1) & \SU_{1,5}/S(U_1U_5) & 11 & 32 &r=2\\ \hline
  E_6^6/\Sp(4) & \SO^0_{5,5}/\SO_5\SO_5 & 26 & 42 &r=6\\
   E_6(\mathbb C)/E_6 & \SO(10,\mathbb C)/\SO(10) & 46 & 78 &r=6\\ \hline
  E_7^7/\SU(8) & E_6^6/\Sp(4) & 43 & 70 &r=7\\
  E_7(\mathbb C)/E_7 & E_6(\mathbb C)/E_6 & 79 & 133 &r=7\\ \hline
  \end{array}
\end{equation*}

\begin{equation*}
\begin{array}{|c|c|c|c|c|}
\hline
  X& Y& \text{srk}(X)&n&\text{Comments}\\ \hline

  E_8^8/\SO(16) & E_7^7/\SU(8) & 71 & 128 &r=8\\
  E_8(\mathbb C)/E_8 & E_7(\mathbb C)/E_7 & 134 & 248 &r=8\\ \hline
  F_4^4/\Sp(3)\times\Sp(1) & \SO^0_{3,4}/\SO_3\SO_4\; \text{or}\;\Sp(3,\mathbb R)/U(3) & 13 & 28 & r=4\\
  E_6^2/\SU(6)\times\Sp(1) & \SU_{3,3}/S(U_3 U_3) & 19 & 40&r=4\\
  E_7^{-5}/\SO(12)\times\Sp(1) & \SO^*(12)/U(6) & 31 & 64 &r=4\\
  E_8^{-24}/E_7\times \Sp(1) & E_7^{-25}/E_6\times U(1) & 55 & 112&r=4\\
  F_4(\mathbb C)/F_4 & \SO(7,\mathbb C)/\SO(7)\;\text{or}\;\Sp(3,\mathbb C)/\Sp(3) & 22 & 52&r=4\\ \hline
  G_2^2/\SO(4) & \mathbb H ^2 & 3 & 8&r=2 \\
  G_2(\mathbb C)/G_2 & \mathbb H ^3 & 4 & 14 &r=2\\ \hline
\end{array}
\end{equation*}

\vskip 5pt

\centerline{\bf Table 1: Splitting rank of irreducible symmetric spaces of non-compact type}

\begin{remark}
In the above table, the symmetric spaces are listed according to their Dynkin diagrams. The groups are listed in the order $A_r$, $B_r$, $C_r$, $D_r$, $(BC)_r$, $E_6$, $E_7$, $E_8$, $F_4$ and $G_2$. Notice that a symmetric space of non-compact type is uniquely determined by its Dynkin diagram together with the multiplicities ($\dim(\mathfrak{p}_\alpha)$) of simple roots.
\end{remark}

\begin{figure}[h]
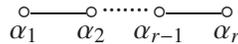

\makebox[1pt][c]{
\beginpicture
\setcoordinatesystem units <.9cm,.9cm> point at -.4 1
\setplotarea x from -.4 to 3.4, y from -1.4 to 1.4
\linethickness=.7pt
\putrule from 0.1 0  to 0.9 0
\putrule from 2.1 0  to 2.9 0
\put {$\cdot$} at 1.2 0
\put {$\cdot$} at 1.3 0
\put {$\cdot$} at 1.4 0
\put {$\cdot$} at 1.5 0
\put {$\cdot$} at 1.6 0
\put {$\cdot$} at 1.7 0
\put {$\cdot$} at 1.8 0
\put {$\circ$} at  0  0
\put {$\circ$} at  1 0
\put {$\circ$} at  2  0
\put {$\circ$} at 3  0
\put {$\alpha_1$} at  0 -.3
\put {$\alpha_2$} at  1 -.3
\put {$\alpha_{r-1}$} at  2 -.3
\put {$\alpha_r$} at  3 -.3
\endpicture
}

\caption{Dynkin diagram of type $A_r$}
\end{figure}

\begin{proof}
We prove the case of $\SL(r+1,\mathbb{R})/\SO(r+1)$ (the first line in Table 1). If $X=\SL(r+1,\mathbb{R})/\SO(r+1)$, the Dynkin diagram is shown in Figure 1, with multiplicities all ones. By the previous discussion, $Y$ is generated by the truncated simple system $\{\alpha_1,...,\hat{\alpha_i},...,\alpha_r\}$ for some $i=1,...,r$, preserving the same multiplicities and configurations. Hence we have $Y=\SL(i,\mathbb{R})/\SO(i)\times \SL(r-i+1,\mathbb{R})/\SO(r-i+1)$, for some $i=1,...,r$. The dimension of $Y$ equals $(i-1)(i+2)/2+(r-i)(r-i+3)/2$, so the codimension of $Y\times \mathbb{R}\subset X$ is $-i^2+(r+1)i$, which attains its minimal codimension $r$ when $i=1,r$. In both cases we have $Y=\SL(r,\mathbb{R})/\SO(r)$, and hence $\text{srk}(X)=\dim(Y\times \mathbb R)=n-r$, where $n=\dim(X)=r(r+3)/2$. The remaining cases are analyzed similarly, and can be found in the Appendix.
\end{proof}

By looking at the table above, we immediately have the following corollaries.

\begin{cor}\label{cor:connect}
Under the same assumptions of Theorem \ref{thm:srk}, $Y$ is also irreducible.
\end{cor}

\begin{cor}\label{cor:si>rank}
If $X$ is an irreducible symmetric space of non-compact type, then $\text{srk}(X)\leq \dim(X)-\rank(X)$, and equality holds if and only if $X=\SL(r+1,\mathbb R)/\SO(r+1)$.
\end{cor}

Before moving to the next section, we will need a little more information about the dimensions of totally geodesic submanifolds $Y\times \mathbb R\subset X$. Here we focus on the cases where $Y$ is of non-compact type, that is, $Y$ is generated by a truncated simple system $\{\alpha_1,...,\hat{\alpha_i},...,\alpha_r\}$ where $\{\alpha_1,...,\alpha_r\}$ is a simple system of $X$. Besides the largest dimension case at $\text{srk}(X)$, we are also curious about the second largest dimension, and we will need to verify that there is enough gap between the two. Such phenomenon of large gap happens to be quite useful when we try to estimate the $k$-th splitting rank. (See below Section \ref{sec:srk^k} and Theorem \ref{thm:brain})
\begin{prop}\label{prop:gap}
Let $X$ be an irreducible symmetric space of non-compact type. Assume $\dim(X)=n$ and $\rank(X)=r\geq 4$. If $Y\times \mathbb R$ and $Y'\times \mathbb R$ are two totally geodesic submanifolds in $X$ where $\dim(Y\times \mathbb R)=\text{srk}(X)$ and $\dim(Y'\times \mathbb R)$ attains the second largest dimension among all truncated simple systems generating $Y'$, then the gaps between the two dimensions ($\dim(Y)-\dim(Y')$) are given in the following table.
\end{prop}

\renewcommand{\arraystretch}{1.5}
\begin{equation*}
\begin{array}{|c|c|c|c|}
\hline
  X & Y' &  \text{Gap}& \text{Comments} \\ \hline
  \SL(r+1,\mathbb{R})/\SO(r+1) & \mathbb H^2\times\SL(r-1,\mathbb{R})/\SO(r-1) &  r-2 &r\geq 4\\
  \SL(r+1,\mathbb{C})/\SU(r+1) & \mathbb H^3\times\SL(r-1,\mathbb{C})/\SU(r-1) & 2r-4 &r\geq 4\\
  \SU^*(2r+2)/\Sp(r+1) & \mathbb H^5\times\SU^*(2r-2)/\Sp(r-1) & 4r-8 &r\geq 4\\ \hline

  \SO^0_{r,r+k}/\SO_r\SO_{r+k}&\mathbb H^2\times \SO^0_{r-2,r-2+k}/\SO_{r-2} \SO_{r-2+k} & 2r+k-5 &r+2k> 7\\
  \SO^0_{4,5}/\SO_4 \SO_5 &\SL(4,\mathbb R)/\SO(4)& 3 &r=4,k=1\\
  \SO^0_{5,6}/\SO_5 \SO_6 &\mathbb H^2\times\SO^0_{3,4}/\SO_3 \SO_4\;\text{or} \;\SL(5,\mathbb R)/\SO(5) & 6 &r=5,k=1\\

  \SO(2r+1,\mathbb{C})/\SO(2r+1) & \mathbb H^3\times\SO(2r-3,\mathbb{C})/\SO(2r-3)  & 4r-8 &r>5\\
  \SO(9,\mathbb{C})/\SO(9) & \SL(4,\mathbb C)/\SU(4)  & 6 &r=4\\
  \SO(11,\mathbb{C})/\SO(11) & \mathbb H^3\times\SO(7,\mathbb{C})/\SO(7) \; \text{or}\; \SL(5,\mathbb C)/\SU(5) & 12 &r=5\\\hline
   \Sp(r,\mathbb{R})/U(r) & \mathbb H^2\times\Sp(r-2,\mathbb{R})/U(r-2) & 2r-4 &r>5  \\
  \Sp(4,\mathbb{R})/U(4) & \SL(4,\mathbb R)/\SO(4) & 3 &r=4  \\
  \Sp(5,\mathbb{R})/U(5) & \mathbb H^2\times\Sp(3,\mathbb{R})/U(3)\; \text{or}\; \SL(5,\mathbb R)/\SO(5) & 6 &r=5  \\ \hline

  \end{array}
\end{equation*}

\begin{equation*}
\begin{array}{|c|c|c|c|}
\hline
   X & Y' &  \text{Gap}& \text{Comments} \\ \hline
  \SU_{r,r}/S(U_r U_r) & \mathbb H^3\times\SU_{r-2,r-2}/S(U_{r-2} U_{r-2})  & 4r-9 &r>6 \\
  \SU_{4,4}/S(U_4 U_4) & \SL(4,\mathbb C)/\SU(4)  & 3 &r=4 \\
  \SU_{5,5}/S(U_5 U_5) &   \SL(5,\mathbb C)/\SU(5) & 8 &r=5 \\
  \SU_{6,6}/S(U_6 U_6) & \mathbb H^3\times\SU_{4,4}/S(U_4 U_4)\;\text{or}\;\SL(6,\mathbb C)/\SU(6)  &15 &r=6 \\
  \Sp(r,\mathbb{C})/\Sp(r) & \mathbb H^3\times\Sp(r-2,\mathbb{C})/\Sp(r-2)  & 4r-8 &r>5 \\
  \Sp(4,\mathbb{C})/\Sp(4) & \SL(4,\mathbb C)/\SU(4)  & 6 &r=4 \\\hline
  \Sp(5,\mathbb{C})/\Sp(5) & \mathbb H^3\times\Sp(3,\mathbb{C})/\Sp(3)\;\text{or}\;\SL(5,\mathbb C)/\SU(5)  & 12 &r=5 \\
  \SO^*(4r)/U(2r) & \mathbb H^4\times\SO^*(4r-8)/U(2r-4) & 8r-19 &r>6\\
  \SO^*(16)/U(8) & \SU^*(8)/\Sp(4) & 3 &r=4\\
  \SO^*(20)/U(10) & \SU^*(10)/\Sp(5) & 12 &r=5\\
  \SO^*(24)/U(12) & \SU^*(12)/\Sp(6) & 25 &r=6\\
  \Sp_{r,r}/\Sp_r \Sp_r & \mathbb H^5\times\Sp_{r-2,r-2}/\Sp_{r-2} \Sp_{r-2} & 8r-17 &r>5\\
  \Sp_{4,4}/\Sp_4\Sp_4 & \SU^*(8)/\Sp(4) & 9 &r=4\\
  \Sp_{5,5}/\Sp_5 \Sp_5 & \SU^*(10)/\Sp(5) & 20 &r=5\\ \hline
  \SO^0_{r,r}/\SO_r \SO_r & \mathbb H^2\times\SO^0_{r-2,r-2}/\SO_{r-2} \SO_{r-2} & 2r-5 & r>7 \\
  \SO^0_{4,4}/\SO_4 \SO_4 & \mathbb H^2\times\mathbb H^2\times\mathbb H^2 & 3 & r=4 \\
  \SO^0_{5,5}/\SO_5\SO_5 & \SL(5,\mathbb R)/\SO(5) & 2 & r=5 \\
  \SO^0_{6,6}/\SO_6 \SO_6 & \SL(6,\mathbb R)/\SO(6) & 5 & r=6 \\
  \SO^0_{7,7}/\SO_7 \SO_7 & \mathbb H^2\times\SO^0_{5,5}/\SO_5 \SO_5\;\text{or}\; \SL(7,\mathbb R)/\SO(7) & 9 & r=7 \\

  \SO(2r,\mathbb{C})/\SO(2r) & \mathbb H^3\times\SO(2r-4,\mathbb{C})/\SO(2r-4) & 4r-10 &r>7 \\
  \SO(8,\mathbb{C})/\SO(8) & \mathbb H^3\times\mathbb H^3\times\mathbb H^3 & 6 &r=4 \\
  \SO(10,\mathbb{C})/\SO(10) & \SL(5,\mathbb C)/\SU(5) & 4 &r=5 \\
  \SO(12,\mathbb{C})/\SO(12) & \SL(6,\mathbb C)/\SU(6) & 10 &r=6 \\
  \SO(14,\mathbb{C})/\SO(14) & \mathbb H^3\times\SO(10,\mathbb{C})/\SO(10)\;\text{or}\;\SL(7,\mathbb C)/\SU(7) & 18 &r=7 \\ \hline
  
\end{array}
\end{equation*}
\begin{equation*}
\begin{array}{|c|c|c|c|}
\hline
  X & Y' &  \text{Gap}& \text{Comments} \\ \hline
  \SU_{r,r+k}/S(U_r U_{r+k}) & \mathbb H^3\times\SU_{r-2,r-2+k}/S(U_{r-2} U_{r-2+k}) & 4r+2k-9 & r+2k>6 \\
  \SU_{4,5}/S(U_4 U_5) & \mathbb H^3\times\SU_{2,3}/S(U_2 U_3)\;\text{or}\;\SL(4,\mathbb C)/\SU(4) & 9 & r=4,k=1 \\
  \Sp_{r,r+k}/\Sp_r \Sp_{r+k} & \mathbb H^5\times\Sp_{r-2,r-2+k}/\Sp_{r-2} \Sp_{r-2+k} & 8r+4k-17 & r\geq 4,k\geq 1\\
  \SO^*(4r+2)/U(2r+1) & \mathbb H^5\times\SO^*(4r-6)/U(2r-3)  & 8r-15 & r>4\\
  \SO^*(18)/U(9) & \SU^*(8)/\Sp(4) & 15 & r=4\\ \hline
  E_6^6/\Sp(4) & \SL(6,\mathbb R)/\SO(6) & 5 &r=6\\
  E_6(\mathbb C)/E_6 & \SL(6,\mathbb C)/\SU(6) & 10 &r=6\\ \hline

  E_7^7/\SU(8) & \SO^0_{6,6}/\SO_6\SO_6  & 6 &r=7\\
  E_7(\mathbb C)/E_7 & \SO(12,\mathbb C)/\SO(12) & 12 &r=7\\ \hline

  E_8^8/\SO(16) & \SO^0_{7,7}/\SO_7\SO_7  & 21 &r=8\\
  E_8(\mathbb C)/E_8 & \SO(14,\mathbb C)/\SO(14) & 42 &r=8\\ \hline

  F_4^4/\Sp(3)\times\Sp(1) & \mathbb H^2\times\SL(3,\mathbb R)/\SO(3)  & 5 &r=4\\
  E_6^2/\SU(6)\times\Sp(1) & \SO^0_{3,5}/\SO_3\SO_5 & 3&r=4\\
  E_7^{-5}/\SO(12)\times\Sp(1) & \SO^0_{3,7}/\SO_3 \SO_7 & 3&r=4 \\
  E_8^{-24}/E_7\times \Sp(1) & \SO^0_{3,11}/\SO_3\SO_{11} & 21&r=4\\
  F_4(\mathbb C)/F_4 & \mathbb H^3\times \SL(3,\mathbb C)/\SU(3)  & 10&r=4\\ \hline
\end{array}
\end{equation*}
\vskip 5pt

\centerline{\bf Table 2: Dimension gap after splitting rank}

\begin{proof}
We prove the case of $\SL(r+1,\mathbb R)/\SO(r+1)$ (the first line of Table 2). As we see in the proof of Theorem \ref{thm:srk}, $Y=\SL(i,\mathbb R)\times \SL(r-i+1,\mathbb R)/\SO(r-i+1)$. And the codimesion of $Y\times \mathbb R\subset X$ is $-i^2+(r+1)i$, which attains its minimal codimension $r$ when $i=1,r$. Now it attains its second minimal value $2r-2$ when $i=2,r-1$ provided $r\geq 3$. In this case, $Y=\mathbb H^2\times \SL(r-1,\mathbb R)/\SO(r-1)$, and the gap is $r-2$. Again, the remaining cases are analyzed similarly in the Appendix.
\end{proof}

\begin{lem}(\text{Gap})\label{lem:gap}
Let $X$ be an irreducible symmetric space of non-compact type. Assume $\dim(X)=n$ and $rank(X)=r\geq 3$. If $Y\times \mathbb R$ is a totally geodesic submanifold whose dimension attains $\text{srk}(X)$, and $Y'\times \mathbb R$ is another totally geodesic submanifold whose dimension is $<\text{srk}(X)$. Then either $Y'$ is irreducible or the gap in dimensions of the two ($\dim(Y)-\dim(Y')$) is at least $r-2$.
\end{lem}

\begin{proof}
For $r=3$, the inequality is automatic. For $r\geq 4$, we check that the inequality follows from Table 2 above for all cases except for $\SO(5,5)/\SO(5)\times \SO(5)$, where the gap between the largest ($\SO(4,4)/\SO(4)\times \SO(4)$) and second largest dimension ($\SL(5,\mathbb R)/\SO(5)$) is $2$. However, the space $\SL(5,\mathbb R)/\SO(5)$ is irreducible. So the gap between $\SO(4,4)/\SO(4)\times \SO(4)$ and any reducible space $Y'$ will be at least $3$. This completes the proof.
\end{proof}

\subsection{The $k$-th Splitting Rank}\label{sec:srk^k}
\begin{defn}\label{defn:srk}
Let $X$ be a rank $r$ symmetric space of non-compact type. For each $k$ ($1\leq k\leq r$), we define the \emph{$k$-th splitting rank} of $X$, denoted $\text{srk}^k(X)$, to be the maximal dimension of a totally geodesic submanifold $Y\subset X$ which splits off an isometric $\mathbb{R}^k$-factor.
\end{defn}
\begin{remark}
In the above notion, the first splitting rank is just our previous notion of splitting rank. We can also see that $\text{srk}^{k+1}(X)\leq \text{srk}^{k}(X)$ for $1\leq k \leq r-1$ and $\text{srk}^r(X)=r$. We abuse notation and set $\text{srk}^0(X)=\dim(X)$.
\end{remark}

\begin{prop}\label{prop:k-max mfld}
Suppose $Y \times \mathbb{R}^k\subset X$ has the maximal dimension, that is, $\dim (Y \times \mathbb{R}^k)=\text{srk}^k(X)$. Then the corresponding Lie triple system $[\mathfrak{p}',[\mathfrak{p}',\mathfrak{p}']]\subset \mathfrak{p}'$ has the form $\mathfrak{p}'=\mathfrak{a}\bigoplus_{\alpha\in \Lambda^+, \alpha(V_k)=0}\mathfrak{p}_\alpha$, where $\mathfrak{a}$ is a choice of maximal abelian subalgebra in $\mathfrak{p}$ that contains the $k$-dimensional Euclidean factor $V_k$.
\end{prop}

\begin{proof}
The proof is the same as Proposition \ref{prop:max}, just replacing $V$ with $V_k$.
\end{proof}

\begin{prop}\label{prop:srk gap}
For $1\leq k \leq r-1$, we have strict inequality $\text{srk}^{k+1}(X)< \text{srk}^{k}(X)$. Therefore, the totally geodesic submanifold $Y\times \mathbb R^k$ that has dimension $\text{srk}^{k}(X)$ does not split off any further $\mathbb R$-factors.
\end{prop}
\begin{proof}
Suppose $Y_{k+1}\times \mathbb{R}^{k+1}\subset X$ has the dimension $\text{srk}^{k+1}(X)$. According to Proposition \ref{prop:k-max mfld}, the tangent space of $Y_{k+1}\times \mathbb{R}^{k+1}$ is identified with $\mathfrak{p}'=\mathfrak{a}\bigoplus_{\alpha\in \Lambda^+, \alpha(V_{k+1})=0}\mathfrak{p}_\alpha$, for some $\mathfrak{a}$ that contains the $\mathbb{R}^{k+1}$-factor $V_{k+1}$. We can choose a pair of root vectors $\pm H_{\alpha_0}$ so that it does not lie in the orthogonal complement $V_{k+1}^\perp$. Let $V_k=V_{k+1}\cap H_{\alpha_0}^\perp$ be a $k$-dimensional subspace of $V_{k+1}$, the Lie triple system $\mathfrak{p}''=\mathfrak{a}\bigoplus_{\alpha\in \Lambda^+, \alpha(V_{k})=0}\mathfrak{p}_\alpha$ strictly contains $\mathfrak{p}'$ since $\alpha_0(V_k)=0$ but $\alpha_0(V_{k+1})\neq0$. Therefore, $\text{srk}^{k+1}(X)=\dim \mathfrak{p}'< \dim \mathfrak{p}''\leq \text{srk}^k(X)$.
\end{proof}

\begin{lem}\label{lem:k-srk estimate}
Let $X$ be an irreducible rank $r$ ($r\geq 2$) symmetric space of non-compact type. Then $\text{srk}^k(X)\leq \text{srk}(X)-2(k-1)$ holds for all $1\leq k< r$.
\end{lem}

\begin{proof}
We show this by induction on the rank of the symmetric space. For $r=2$, the only possible value for $k$ is $k=1$, and the inequality holds immediately. Suppose we have the inequality for all such irreducible symmetric spaces of rank $l$ ($l\geq 2$), assuming $\text{rank}(X)=l+1$, we want to show $\text{srk}^k(X)\leq \text{srk}(X)-2(k-1)$ for all $1\leq k< l+1$. Notice when $k=1$, the inequality is trivially true, so we may assume $k\geq 2$.

Let $\text{srk}^k(X)=\dim(Y_k\times \mathbb R^k)$, where $Y_k$ is described as in Proposition \ref{prop:k-max mfld}. Let $V_k$ denote the $\mathbb R^k$-factor. We inductively define $V_i$ so that it is an $i$ dimensional Euclidean subspace of $V_{i+1}$ and $\ker(V_{i+1})\cap \Lambda\subsetneqq \ker(V_i)\cap \Lambda$. This will give rise to an extending chain of Lie triples $\mathfrak{p}_k\subset ...\subset \mathfrak{p}_1\subset\mathfrak{p}$, corresponding to a totally geodesic chain $Y_k\times \mathbb R^k\subset ...\subset Y_1\times \mathbb R\subset X$, such that for each $1\leq i\leq k$, $\mathfrak p_i=\mathfrak{a}\bigoplus_{\alpha\in \Lambda^+, \alpha(V_i)=0}\mathfrak{p}_\alpha$. The choice of $V_i$ implies that $\dim(\mathfrak p_{i+1})< \dim(\mathfrak p_i)$ and therefore $Y_i$ does not split off an $\mathbb R$-factor, for all $i$. Besides, since $Y_k\times \mathbb R^k\subset Y_1\times \mathbb R$ have a common $\mathbb R$-factor $V_1$, $Y_k\times \mathbb R^{k-1}$ is totally geodesic in $Y_1$.

Now if $Y_1$ is irreducible, by the induction hypothesis, we have $\dim(Y_k\times \mathbb R^{k-1})\leq \text{srk}^{k-1}(Y_1)\leq \text{srk}(Y_1)-2(k-2)$. According to Corollary \ref{cor:si>rank}, we have $\text{srk}(Y_1)\leq \dim(Y_1)-l\leq \dim(Y_1)-2$. Hence combining the two inequalities, we conclude $\dim(Y_k\times \mathbb R^k)\leq\text{srk}(Y_1)-2(k-2)+1\leq \dim(Y_1)-2-2(k-2)+1=\dim(Y_1\times \mathbb R)-2(k-1)\leq\text{srk}(X)-2(k-1)$.

If $Y_1$ is reducible, then $Y_1\times \mathbb R$ can not have the dimension equal to the splitting rank of $X$ by Corollary \ref{cor:connect}. So Lemma \ref{lem:gap} implies that $\dim(Y_1\times \mathbb R)\leq \text{srk}(X)-(l+1-2)$. The increasing chain $\mathfrak{p}_k\subset ...\subset \mathfrak{p}_1\subset\mathfrak{p}$ gives the inequality $\dim(\mathfrak{p}_k)\leq \dim(\mathfrak{p}_1)-(k-1)$. Notice $\mathfrak{p}_i$ is identified with the tangent space of $Y_i\times \mathbb R^i$, so we can then estimate $\dim(Y_k\times \mathbb R^k)=\dim(\mathfrak{p}_k)\leq \text{srk}(X)-(l+1-2)-(k-1)\leq \text{srk}(X)-2(k-1)$.

We have shown in both cases that the inequality holds for all rank $l+1$ irreducible symmetric spaces. This completes the induction argument and hence the proof of this lemma.
\end{proof}

\begin{cor}\label{cor:k-srk estimate}
Let $X$ be an irreducible rank $r$ ($r\geq 1$) symmetric space of non-compact type, excluding $\SL(3,\mathbb R)/\SO(3)$, $\Sp(2,\mathbb R)/U(2)$, $G_2^2/\SO(4)$ and $\SL(4,\mathbb R)/\SO(4)$. Then $\text{srk}^k(X)\leq \text{srk}(X)-2(k-1)$ holds for all $1\leq k\leq r$.
\end{cor}
\begin{proof}
Notice the inequality automatically holds in rank one cases, and in view of Lemma \ref{lem:k-srk estimate} we only need to consider the case when $k=r$. If $k=r$, the inequality is equivalent to $3r-2\leq \text{srk}(X)$. We know that the dimension ($=n$) of an irreducible symmetric space grows roughly quadratically in its rank ($=r$), and actually we can check that $n\geq 3r$ whenever $r\geq 3$. This together with Corollary \ref{cor:connect} implies that $\text{srk}(X)=\dim(Y\times \mathbb R)\geq 3(r-1)+1=3r-2$ whenever $r\geq 4$. This proves the corollary in all cases where $\rank\geq 4$. When $r=2$, it is equivalent to show $\text{srk}(X)\geq 4$, and according to Table 1, this excludes $\SL(3,\mathbb R)/\SO(3)$, $\Sp(2,\mathbb R)/U(2)$ and $G_2^2/\SO(4)$. When $r=3$, it is equivalent to show $\text{srk}(X)\geq 7$, and by Table 1 this only excludes $\SL(4,\mathbb R)/\SO(4)$.
\end{proof}

\subsection{Reducible Symmetric Spaces}\label{sec:reducible}
We intend to generalize Corollary \ref{cor:k-srk estimate} to all higher rank symmetric spaces of non-compact type. Below is the key lemma that characterizes certain $\mathbb R$-split totally geodesic submanifolds in reducible symmetric spaces.

\begin{lem}\label{lem:splitting}
Let $X$ be a symmetric space of non-compact type, and $Z=Y\times \mathbb R^k$ a totally geodesic subspace in $X$ that has dimension equal to the $k$-th splitting rank of $X$. If $X$ splits as a product of $X_1$ and $X_2$, then $Z$ also splits as a product of $Z_1$ and $Z_2$, where $Z_i=Y_i\times \mathbb R^{k_i}$ is totally geodesic in $X_i$, for $i=1,2$ and some $k_i\geq 0$ satisfying $k_1+k_2=k$.
\end{lem}
\begin{proof}
We write $X=G/K$ and fix a basepoint $x\in X$. We form the Cartan decomposition $\mathfrak{g}=\mathfrak{k}+\mathfrak{p}$, where $\mathfrak{p}$ can be identified with the tangent space $T_xX$. We denote $V_k\subset \mathfrak{p}$ the $\mathbb R^k$ factor of $Z$, and extend $V_k$ to a maximal abelian subalgebra $\mathfrak{a}\subset \mathfrak{p}$. Since $X$ splits as a product of $X_1$ and $X_2$, we can write $\mathfrak{a}=\mathfrak{a}_1\oplus \mathfrak{a}_2$, and also the set of roots $\Lambda$ of $X$ decomposes as $\Lambda_1\cup \Lambda_2$, where $\Lambda_i$ is the set of roots belonging to $X_i$ with respect to $\mathfrak{a}_i$. By Proposition \ref{prop:k-max mfld}, $Z$ has the Lie triple system $\mathfrak{p}'=\mathfrak{a}\bigoplus_{\alpha\in \Lambda^+, \alpha(V_k)=0}\mathfrak{p}_\alpha$. Let $\mathfrak{p}_i'=\mathfrak{a}_i\bigoplus_{\alpha\in \Lambda_i^+, \alpha(V_k)=0}\mathfrak{p}_\alpha$, we have $\mathfrak{p}'=\mathfrak{p}_1'\oplus \mathfrak{p}_2'$. Notice $\mathfrak{p}_i'\subset \mathfrak{p}_i$ is a Lie triple system. Indeed, $\mathfrak{p}_i'=\mathfrak{a}_i\bigoplus_{\alpha\in \Lambda_i^+, \alpha(V_{k,i})=0}\mathfrak{p}_\alpha$ where $V_{k,i}$ is the orthogonal projection of $V_k$ to $\mathfrak{a}_i$. This implies that $V_{k,i}$ is the Euclidean factor of $\mathfrak{p}_i$, therefore $V_k$ splits as $V_{k,1}\oplus V_{k,2}$, which completes the proof.
\end{proof}

\begin{cor}\label{cor:reducible k-srk}
Let $X$ be a rank $r$ symmetric space of non-compact type. Assume $X=X_1\times X_2$, where $\text{rank}(X_i)=r_i$ for $i=1,2$. Then $$\text{srk}^k(X)=\text{Max}\{\text{srk}^{j_1}(X_1)+\text{srk}^{j_2}(X_2):{0\leq j_1\leq r_1,0\leq j_2\leq r_2}\;\text{and}\;j_1+j_2=k\}.$$
\end{cor}
\begin{remark}
This is a direct consequence of Lemma \ref{lem:splitting}. In the corollary, recall that by definition $\text{srk}^0(X)=\dim(X)$. As a result, $\text{srk}(X_1\times X_2)=\text{Max}\{\text{srk}(X_1)+\dim(X_2),\text{srk}(X_2)+\dim(X_1)\}$. Furthermore, if we define \emph{the k-th splitting index of $X$} to be $\text{si}^k(X):=n-\text{srk}^k(X)$, 
 then Corollary \ref{cor:reducible k-srk} simply says $$\text{si}^k(X_1\times X_2)=\text{Min}\{\text{si}^{j_1}(X_1)+\text{si}^{j_2}(X_2):0\leq j_1\leq r_1,0\leq j_2\leq r_2\;\text{and}\;j_1+j_2=k\}$$
and that $\text{si}(X_1\times X_2)=\text{Min}\{\text{si}(X_1),\text{si}(X_2)\}$, where $\text{si}(X)$ denotes the first splitting index of $X$.
\end{remark}

\begin{thm}\label{thm:brain}
Let $X$ be a rank $r$ symmetric space of non-compact type. Assume $X$ has no direct factors isometric to $\mathbb H^2$, $\SL(3,\mathbb R)/\SO(3)$, $\Sp(2,\mathbb R)/U(2)$, $G_2^2/\SO(4)$ or $\SL(4,\mathbb R)/\SO(4)$. Then $\text{srk}^k(X)\leq \text{srk}(X)-2(k-1)$ for all $1\leq k\leq r$.
\end{thm}
\begin{proof}
We write $X$ as a product of irreducible symmetric spaces $X_1\times...\times X_s$. Using the notion of splitting index described in the previous remark, the inequality is equivalent to $\text{si}^k(X)\geq \text{si}(X)+2(k-1)$. By repeatedly applying Corollary \ref{cor:reducible k-srk}, we can assume $\text{si}^k(X)=\sum_{l=1}^s\text{si}^{j_l}(X_l)$ for some $j_l$ satisfying $0\leq j_l\leq r_l$ and $\sum_{l=1}^s j_l=k$ where $r_l$ is the rank of $X_l$. For each $j_l>0$, we have $\text{si}^{j_l}(X_l)\geq \text{si}(X_l)+2(j_l-1)$ by Corollary \ref{cor:k-srk estimate}. Notice $\text{si}^0(X)=0$ and $\text{si}^{j_l}(X_l)$ does not contribute to the summation if $j_l=0$. We can  further estimate

$$\text{si}^k(X)=\sum_{1\leq l\leq s,j_l>0}\text{si}^{j_l}(X_l)\geq \sum_{1\leq l\leq s,j_l>0} [\text{si}(X_l)+2(j_l-1)]= 2k+ \sum_{1\leq l\leq s,j_l>0}(\text{si}(X_l)-2).$$

As a consequence of Corollary \ref{cor:si>rank}, we have $\text{si}(X_l)\geq 2$ as we assume no $\mathbb H^2$-factors. Now we apply the inequality $\text{si}(X_l)\geq \text{Min}_{1\leq l\leq s}\text{si}(X_l)=\text{si}(X)$ to one of the $l$ in the summation, and apply $\text{si}(X_l)\geq 2$ to the rest of the $l$. We finally obtain $\text{si}^k(X)\geq \text{si}(X)+2k-2$, which completes the proof.
\end{proof}

\section{Application to Bounded Cohomology}\label{sec:BC}
In this section, we generalize the method of \cite{LW}, and show the surjectivity of comparison maps in a slightly larger range ($\geq \text{srk}+2$). The approach is quite similar: Theorem \ref{thm:brain} generalizes \cite[Lemma 4.6]{LW}, allowing us to improve its main theorem to the following.

\begin{thm}\label{thm:main}
Let $X=G/K$ be an $n$-dimensional symmetric space of non-compact type of rank $r\geq 2$, and $\Gamma$ a cocompact torsion-free lattice in $G$. Assume $X$ has no direct factors of $\mathbb H^2$, $\SL(3,\mathbb R)/\SO(3)$, $\Sp(2,\mathbb R)/U(2)$, $G_2^2/\SO(4)$ and $\SL(4,\mathbb R)/\SO(4)$, then the comparison maps
$\eta:H^{*}_{c,b}(G,\mathbb{R})\rightarrow H^{*}_c(G,\mathbb{R})$ and $\eta':H^{*}_{b}(\Gamma,\mathbb{R})\rightarrow
H^{*}(\Gamma,\mathbb{R})$ are both surjective in all degrees $* \geq \text{srk}(X)+2$.
\end{thm}

We summarize the approach of \cite{LW} in the following steps.

\textbf{Step 1}: Notice that surjectivity of $\eta'$ implies surjectivity of $\eta$. In order to show surjectivity of $\eta'$ in degree $k$, it is equivalent to assign each cohomology class $[f]\in H^k(\Gamma,\mathbb R)$ a bounded representative. By the explicit isomorphism $H_{dR}^k(X/\Gamma, \mathbb{R})\simeq H_{sing}^k(X/\Gamma,\mathbb{R})\simeq H^k(\Gamma,\mathbb{R})$, we can view $f:\Gamma^k\rightarrow \mathbb R$ as a function that integrates a $k$-form $f_\omega$ over a $k$-simplex generated by a $k$-tuple in $\Gamma$. We now replace it (within the same cohomology class) by a function that integrates the same form $f_\omega$ over a ``barycentrically straightened'' $k$-simplex, and claim it is bounded when $k$ is in certain degrees. This produces the bounded representative. (See \cite[section 5.1, 5.2]{LW} for more details.)

\textbf{Step 2}: For each $1\leq k\leq n$, the barycentric $k$-straightening is a map $st_k: C^k(X)\rightarrow C^k(X)$, where $X=G/K$ is a symmetric space of non-compact type. It gives a $C^1$-chain map and is chain homotopic to the identity. It is also $G$-equivariant hence preserves the $\Gamma$-action. Moreover, if one can in addition show that the straightened $k$-simplices have uniformly bounded Jacobian, then the function $f_\omega$ is also bounded in degree $k$. As a result the surjectivity of the comparison map is obtained in degree $k$. (See \cite[section 2.3]{LW} for more details.) Notice in \cite{LW}, they showed uniformly bounded Jacobian for irreducible symmetric spaces excluding $\SL(3,\mathbb R)/\SO(3)$ and $\SL(4,\mathbb R)/\SO(4)$ when $k\geq n-r+2$. This is generalized in this paper, where we show uniformly bounded Jacobian in degrees $k\geq \text{srk}(X)+2$, for all symmetric spaces satisfying the condition of Theorem \ref{thm:brain}.

\textbf{Step 3}: By the computation in \cite[section 2.3]{LW}, the Jacobian of straightened $k$-simplex is bounded above (up to a multiplicative constant) by the quotient $\det(Q_1|_S)^{1/2}/\det(Q_2|_S)$, where $S$ is a $k$-dimensional subspace in a tangent space $T_xX$, and $Q_1$, $Q_2$ are two positive semidefinite quadratic forms ($Q_2$ is actually positive definite) defined by the following:
$$Q_1(v,v)=\int_{\partial_F X}dB^2_{(x,\theta)}(v)d\mu(\theta),$$
$$Q_2(v,v)=\int_{\partial_F X}DdB_{(x,\theta)}(v,v)d\mu(\theta).$$
In the above expression, $B$ is the Busemann function on $X$ based at some fixed point, and $\mu$ is a probability measure fully supported on the Furstenberg boundary $\partial_F X$. Therefore if we can bound $\det(Q_1|_S)^{1/2}/\det(Q_2|S)$ by some constant $C$ that only depends on $X$ (independent of the choices $x\in X$ and $S\subset T_xX$), then we are able to control the Jacobian in degree $k=\dim(S)$. (See \cite[section 2.2, 2.3]{LW} for more details.)

\textbf{Step 4}: In order to show the ratio $\det(Q_1|_S)^{1/2}/\det(Q_2|_S)$ is uniformly bounded, we need an eigenvalue matching property. If $S$ is in top dimension $n$, then the ratio is just $\det(Q_1)^{1/2}/\det(Q_2)$. Following the approach from \cite{CF1}, \cite{CF2}, it suffices to give for each small eigenvalues of $Q_2$, two comparably small eigenvalues of $Q_1$ to cancel with. Generalizing this argument, if we were able to find $r-2$ additional small eigenvalues of $Q_1$ to cancel with the smallest eigenvalues of $Q_2$, then we can restrict the quadratic forms to a subspace $S$ of dimension $k$ (where $k\geq n-r+2$), and the ratio of the determinants $\det(Q_1|_S)^{1/2}/\det(Q_2|_S)$ remains uniformly bounded. This is implied by a weak eigenvalue matching theorem. (See \cite[Theorem 3.3]{LW} and originally \cite{CF2}.) Actually, there are at most $r$ many eigenvectors of $Q_2$ with small eigenvalues, and they almost lie in a tangent space of a flat (have small angle with a flat). For each such eigenvector $v_i$, we can always find two unit vectors $v_i'$, $v_i''$ so that the $Q_1$ values on these two vectors are bounded above by $Q_2(v_i,v_i)$, up to a uniform multiplicative constant (we say in this case $v_i'$ and $v_i''$ cancels with $v_i$). Moreover, for the smallest eigenvector $v_1$, we are able to find $r-2$ additional unit vectors $v_1^{(3)},...,v_1^{(r)}$ to cancel with, and the collection of all the pairs of vectors, together with the $r-2$ additional vectors, almost form an orthonormal frame. Then by the Gram-Schmidt process and standard linear algebra, we can find the eigenvalue match. (See \cite[section 3.3]{LW} for more details.)

\textbf{Step 5}: Finally, the weak eigenvalue matching theorem described in Step 4 can be further reduced to a combinatorial problem. For a single vector $v$ that lies in a tangent space $\mathfrak{a}$ of a flat, we denote by $v^*$ the most singular vector in a fixed small neighborhood of $v$. Then any vector that lies in $Q_v=\bigoplus_{\alpha\in \Lambda^+, \alpha(v^*)\neq0}\mathfrak{p}_\alpha$ will be able to cancel $v$. Now for each root $\alpha$, we pick an orthonormal frame $\{b_{\alpha_i}\}$ for $\mathfrak{p}_\alpha$, and we collect them into the set $B=\{b_i\}_{i=1}^{n-r}$. The idea is to find among the set $B$, $3r-2$ distinct vectors $v_1',v_1'',...,v_1^{(r)},v_2',v_2'',...,v_r',v_r''$ to cancel a given almost orthonormal frame $\{v_1,...v_r\}$ in $\mathfrak{a}$ (hence $v_1^*,...v_r^*$ are distinct). Notice this is a generalization of the classic Hall's Marriage Problem, and in order to solve this, it is sufficient to solve a cardinality inequality: for any subcollection of vectors $\{v_{i_1},...,v_{i_k}\}$, the number of vectors in $B$ that belongs to $\bigcup_{j=1}^k Q_{v_{i_j}}$ is at least $2k+r-2$ (which is showed in \cite[Lemma 4.6]{LW}). This solves the eigenvalue matching in the special case where the small eigenvectors of $Q_2$ all lie in a same $\mathfrak{a}$. For the general case where the small eigenvectors have small angles to a flat, a similar argument is used, and we refer the readers to \cite[Section 4]{LW} for more details.\\

\emph{Proof of Theorem \ref{thm:main}}: The proof is similar to that of \cite[Main Theorem]{LW}. Notice that Step 1-3 goes through unchanged in our case, and in Step 4-5 we only need to find $n-\text{srk}(X)-2$ additional (instead of $r-2$) small eigenvalues of $Q_1$ to cancel with the smallest eigenvalue of $Q_2$. Showing this requires a similar weak eigenvalue matching theorem as \cite[Theorem 3.3]{LW}, where the existing $(2k+r-2)$-frame $\{v_1',v_1'',...,v_1^{(r)},v_2',v_2'',...,v_k',v_k''\}$ is now replaced by a $(2k+n-\text{srk}(X)-2)$-frame $\{v_1',v_1'',...,v_1^{(n-\text{srk}(X))},v_2',v_2'',...,v_k',v_k''\}$ and the corresponding angle inequalities are satisfied. This is further reduced to a similar Hall's Marriage type combinatorial problem, and the corresponding cardinality estimate is ensured by the following modification of \cite[Lemma 4.6]{LW}:

\begin{lem}\label{lem:revise}
Let $X=G/K$ be a rank $r\geq 2$ symmetric space of non-compact type, without direct factors isometric to $\mathbb H^2$, $\SL(3,\mathbb R)/\SO(3)$, $\Sp(2,\mathbb R)/U(2)$, $G_2^2/\SO(4)$, or $\SL(4,\mathbb R)/\SO(4)$. Fix a maximal abelian subalgebra $\mathfrak{a}\subset \mathfrak{p}$. Assume $\{v_1^*,...,v_r^*\}$ spans $\mathfrak{a}$, and let $Q_i=\bigoplus_{\alpha\in \Lambda^+, \alpha(v_i^*)\neq0}\mathfrak{p}_\alpha$. Then for any subcollection of vectors $\{v_{i_1}^*,...,v_{i_k}^*\}$, we have $\dim(Q_{i_1}+...+Q_{i_k})\geq (2k+n-\text{srk}(X)-2)$.
\end{lem}

\begin{proof}
Notice $Q_{i_1}+...+Q_{i_k}=\bigoplus_{\alpha\in \Lambda^+, \alpha(V_k)\neq0}\mathfrak{p}_\alpha$ where $V_k$ is the span of $v_{i_1}^*,...,v_{i_k}^*$. Its orthogonal complement in $\mathfrak{p}$ is $\mathfrak{a}\bigoplus_{\alpha\in \Lambda^+, \alpha(V_k)=0}\mathfrak{p}_\alpha$, which has dimension at most $\text{srk}^k(X)$, hence according to Theorem \ref{thm:brain} is bounded above by $\text{srk}(X)-2k+2$. Therefore, $\dim(Q_{i_1}+...+Q_{i_k})\geq (2k+n-\text{srk}(X)-2)$. This completes Lemma \ref{lem:revise} and hence Theorem \ref{thm:main}.
\end{proof}

\begin{remark}
Notice that the Main Theorem in \cite{LW} required the symmetric space to be irreducible. But the proof of Lemma 4.6 was the only step that used irreduciblity. The rest of the proof remains valid for reducible symmetric spaces. Thus replacing \cite[Lemma 4.6]{LW} by our Lemma \ref{lem:revise}, the actual proof goes through even in the reducible case.
\end{remark}

\section{Appendix}
In this section, we finish the proof of Theorem \ref{thm:srk} and Proposition \ref{prop:gap}. We combine the two proofs as they are both a case by case argument. And we notice the case of $\SL(r+1,\mathbb{R})/\SO(r+1)$ has already been proved in the context.

\emph{Case of $\SL(r+1,\mathbb{C})/\SU(r+1)$}: the Dynkin diagram is of type $A_r$ and is shown in Figure 1, with multiplicities $2$ for all simple roots. Since $Y$ is generated by the truncated simple system $\{\alpha_1,...,\hat{\alpha_i},...,\alpha_r\}$ for some $i=1,...,r$, preserving the same multiplicities and configurations, we have $Y=\SL(i,\mathbb{C})/\SU(i)\times \SL(r-i+1,\mathbb{C})/\SU(r-i+1)$, for some $i=1,...,r$. The dimension of $Y$ equals $(i-1)(i+1)+(r-i)(r-i+2)$, so the codimension of $Y\times \mathbb{R}\subset X$ is $-2i^2+2(r+1)i$, which attains its minimal codimension $2r$ when $i=1,r$. In both cases we have $Y=\SL(r,\mathbb{C})/\SO(r)$, and hence $\text{srk}(X)=\dim(Y\times \mathbb R)=n-2r$, where $n=\dim(X)=r(r+2)$. The codimsion of $Y\times \mathbb{R}\subset X$ attains its second minimal value when $i=2, r-1$ provided $r\geq 3$. In this case, $Y=\mathbb H^3\times \SL(r-1,\mathbb{C})/\SU(r-1)$ and the codimension of $Y\times \mathbb{R}$ is $4r-4$, so the gap is $2r-4$.

\emph{Case of $\SU^*(2r+2)/\Sp(r+1)$}: the Dynkin diagram is of type $A_r$ and is shown in Figure 1, with multiplicities $4$ for all simple roots. Hence $Y=\SU^*(2i)/\Sp(i)\times \SU^*(2r+2-2i)/\Sp(r+1-k)$ for some $i=1,...,r$. The codimension of $Y\times \mathbb{R}\subset X$ is $-4i^2+4(r+1)i$, which attains minimal when $i=1,r$. In both cases $Y=\SU^*(2r)/\Sp(r)$, and $\text{srk}(X)=n-4r$. The codimension attains its second minimal value when $i=2,r-1$ provided $r\geq 3$, in which case $Y=\mathbb H^5\times\SU^*(2r-2)/\Sp(r-1)$. The codimension of $Y\times \mathbb{R}$ is $8r-8$, so the gap is $4r-8$.

\emph{Case of $E_6^{-26}/F_4$}: the Dynkin diagram is of type $A_2$ and is shown in Figure 1 where $r=2$, with multiplicities $8$ for both simple roots. Hence $Y$ can only be $\mathbb H^9$ so that $\text{srk}(X)=10$. Notice that $E_6^{-26}/F_4$ is of rank two and it does not satisfy the condition of Proposition \ref{prop:gap}.

\begin{figure}
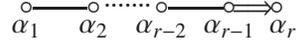

\makebox[1pt][c]{
\beginpicture
\setcoordinatesystem units <.9cm,.9cm> point at -.4 1
\setplotarea x from -.4 to 4.4, y from -1.4 to 1.4
\linethickness=.7pt
\putrule from 0.1 0  to 0.9 0
\putrule from 2.1 0  to 2.9 0
\put {$\cdot$} at 1.2 0
\put {$\cdot$} at 1.3 0
\put {$\cdot$} at 1.4 0
\put {$\cdot$} at 1.5 0
\put {$\cdot$} at 1.6 0
\put {$\cdot$} at 1.7 0
\put {$\cdot$} at 1.8 0
\put {$\circ$} at  0  0
\put {$\circ$} at  1 0
\put {$\circ$} at  2  0
\put {$\circ$} at 3  0
\put {$\circ$} at 3.7  0
\put {$\Longrightarrow$} at 3.35 0
\put {$\alpha_1$} at  0 -.3
\put {$\alpha_2$} at  1 -.3
\put {$\alpha_{r-2}$} at  2 -.3
\put {$\alpha_{r-1}$} at  3 -.3
\put {$\alpha_r$} at 3.8 -.3
\endpicture
}
\caption{Dynkin diagram of type $B_r$}
\end{figure}

\emph{Case of $\SO_0(r,r+k)/\SO(r)\times \SO(r+k)$}: the Dynkin diagram is of type $B_r$ and is shown in Figure 2, with ordered multiplicities $1,1,...,1,k$. If we remove $\alpha_i$, the remaining diagram (with multiplicity information) will represent $Y_i=\SL(i,\mathbb R)/\SO(i)\times [\SO_0(r-i,r-i+k)/\SO(r-i)\times \SO(r-i+k)]$ (Notice $\SL(1,\mathbb R)/\SO(1)$ and  $\SO(0,k)_0/\SO(0)\times \SO(k)$ are just a point by abuse of notation). Thus we can compute that $Y_i\times \mathbb R$ has codimension $-3i^2/2+(4r+2k-1)i/2$ in $X$. It attains a minimum when $i=1$ provided $r+2k>4$, and so $Y=\SO_0(r-1,r-1+k)/\SO(r-1)\times \SO(r-1+k)$, $\text{srk}(X)=n-(2r+k-2)$. The only space that satisfies $r+2k\leq 4$ is $\SO_0(2,3)/\SO(2)\times \SO(3)$, and it has splitting rank $3$, corresponding to $\mathbb H^2\times \mathbb R$. This agrees with the general formula hence can be absorbed into it. Now the codimension attains its second minimum when $i=2$ provided $r+2k>7$, and so $Y'=\mathbb H^2\times\SO_0(r-2,r-2+k)/\SO(r-2)\times \SO(r-2+k)$, $\dim(Y'\times \mathbb R)=n-(4r+2k-7)$. Hence the gap is $2r+k-5$. As we focus on $r\geq 4$ in Proposition \ref{prop:gap}, the spaces that are excluded by $r+2k>7$ is $\SO_0(4,5)/\SO(4)\times \SO(5)$ and $\SO_0(5,6)/\SO(5)\times \SO(6)$. If $X=\SO_0(4,5)/\SO(4)\times \SO(5)$ ($\dim(X)=20$), then $Y$ is $\SO_0(3,4)/\SO(3)\times \SO(4)$ (dimension $12$), and $Y'$ is $\SL(4,\mathbb R)/\SO(4)$ (dimension $9$), so the gap is $3$. If $X=\SO_0(5,6)/\SO(5)\times \SO(6)$, then $Y=\SO_0(4,5)/\SO(4)\times \SO(5)$, and $Y'$ is either $\mathbb H^2\times \SO(3,4)/\SO(3)\times\SO(4)$ or $\SL(5,\mathbb R)/\SO(5)$ with dimension $14$. So the gap is $6$.

\emph{Case of $\SO(2r+1,\mathbb{C})/\SO(2r+1)$}: the Dynkin diagram is of type $B_r$ and is shown in Figure 2, with multiplicities $2$ for all simple roots. If we remove $\alpha_i$, the remaining diagram will represent $Y_i=\SL(i,\mathbb C)/\SU(i)\times \SO(2r-2i+1,\mathbb{C})/\SO(2r-2i+1)$ (notice $\SL(1,\mathbb C)/\SU(1)$ and $\SO(1,\mathbb C)/\SO(1)$ are just a point by abuse of notation). We compute that the codimension of $Y_i\times \mathbb R$ in $X$ is $-3i^2+(4r+1)i$. It has minimal value $4r-2$ when $i=1$ (when $r=2$ it takes minimal value on both $i=1,2$, but they both represent the same subspace $\mathbb H^3$). Hence the splitting rank is $n-(4r-2)$, corresponding to the subspace $Y=\SO(2r-1,\mathbb C)/\SO(2r-1)$. Now the codimension takes the second minimal value $8r-10$ when $i=2$, provided $r>5$. In this case, $Y'=\mathbb H^3\times \SO(2r-3,\mathbb C)/\SO(2r-3)$ and the gap is $4r-8$. If $r=4$, then $X=\SO(9,\mathbb C)/\SO(9)$ and $Y=\SO(7,\mathbb C)/\SO(7)$. The codimension takes its second minimal value $20$ when $Y'=\SL(4,\mathbb C)/\SU(4)$, hence the gap is $6$. If $r=5$, then $X=\SO(11,\mathbb C)/\SO(11)$ and $Y=\SO(9,\mathbb C)/\SO(9)$. The codimension takes its second minimal value $30$ when $Y'$ is $\mathbb H^3\times \SO(7,\mathbb C)/\SO(7)$ or $\SL(5,\mathbb C)/\SU(5)$, hence the gap is $12$.

\begin{figure}
\makebox[1pt][c]{
\beginpicture
\setcoordinatesystem units <.9cm,.9cm> point at -.4 1
\setplotarea x from -.4 to 4.4, y from -1.4 to 1.4
\linethickness=.7pt
\putrule from 0.1 0  to 0.9 0
\putrule from 2.1 0  to 2.9 0
\put {$\cdot$} at 1.2 0
\put {$\cdot$} at 1.3 0
\put {$\cdot$} at 1.4 0
\put {$\cdot$} at 1.5 0
\put {$\cdot$} at 1.6 0
\put {$\cdot$} at 1.7 0
\put {$\cdot$} at 1.8 0
\put {$\circ$} at  0  0
\put {$\circ$} at  1 0
\put {$\circ$} at  2  0
\put {$\circ$} at 3  0
\put {$\circ$} at 3.7  0
\put {$\Longleftarrow$} at 3.35 0
\put {$\alpha_1$} at  0 -.3
\put {$\alpha_2$} at  1 -.3
\put {$\alpha_{r-2}$} at  2 -.3
\put {$\alpha_{r-1}$} at  3 -.3
\put {$\alpha_r$} at 3.8 -.3
\endpicture
}
\caption{Dynkin diagram of type $C_r$}
\end{figure}

\emph{Case of $\Sp(r,\mathbb{R})/U(r)$}: the Dynkin diagram is of type $C_r$ and is shown in Figure 3, with multiplicities $1$ for all simple roots. If we remove $\alpha_i$, the remaining diagram will represent $Y_i=\SL(i,\mathbb R)/\SO(i)\times \Sp(r-i,\mathbb{R})/U(r-i))$ (notice $\SL(1,\mathbb R)/\SO(1)$ and $\Sp(0,\mathbb{R})/U(0)$ are just a point by abuse of notation). We compute that the codimension of $Y_i\times \mathbb R$ in $X$ is $-3i^2/2+(4r+1)i/2$. It has minimal value $2r-1$ when $i=1$ (when $r=2$ it takes minimal value on both $i=1,2$, but they both represent the same space $\mathbb H^2$). Hence the splitting rank is $n-(2r-1)$ corresponding to the space $Y=\Sp(r-1,\mathbb{R})/U(r-1)$. Now the codimension takes the second minimal value $4r-5$ when $i=2$, provided $r>5$. In this case, $Y'=\mathbb H^2\times \Sp(r-2,\mathbb{R})/U(r-2)$ and the gap is $2r-4$. If $r=4$, then $X=\Sp(4,\mathbb{R})/U(4)$ and $Y=\Sp(3,\mathbb{R})/U(3)$. The codimension takes its second minimal value $10$ when $Y'=\SL(4,\mathbb R)/\SO(4)$, hence the gap is $3$. If $r=5$, then $X=\Sp(5,\mathbb{R})/U(5)$ and $Y=\Sp(4,\mathbb{R})/U(4)$. The codimension takes its second minimal value $15$ when $Y'$ is $\mathbb H^2\times \Sp(3,\mathbb{R})/U(3)$ or $\SL(5,\mathbb R)/\SO(5)$, hence the gap is $6$.

\emph{Case of $\SU(r,r)/S(U(r)\times U(r))$}: the Dynkin diagram is of type $C_r$ and is shown in Figure 3, with ordered multiplicities $2,2,...,2,1$. If we remove $\alpha_i$, the remaining diagram will represent $Y_i=\SL(i,\mathbb C)/\SU(i)\times \SU(r-i,r-i)/S(U(r-i)\times U(r-i))$ (notice $\SL(1,\mathbb C)/\SU(1)$ and $\SU(0,0)/S(U(0)\times U(0))$ are just a point by abuse of notation). We compute that the codimension of $Y_i\times \mathbb R$ in $X$ is $-3i^2+4ri$. It has minimal value $4r-3$ when $i=1$, provided $r>3$. If $r=2$, $X$ is isomorphic to $\SO_0(2,4)/\SO(2)\times\SO(4)$, which has been solved previously. If $r=3$, the codimension is minimal for both $i=1,3$, which corresponds to $\SO_0(2,4)/\SO(2)\times\SO(4)$ and $\SL(3,\mathbb C)/\SU(3)$. Hence the splitting rank is $n-(4r-3)$ corresponding to the space $Y=\Sp(r-1,\mathbb{R})/U(r-1)$ (except for the case $r=3$ where there are two subspaces). Now the codimension takes the second minimal value $8r-12$ when $i=2$, provided $r>6$. In this case, $Y'=\mathbb H^3\times \SU(r-2,r-2)/S(U(r-2)\times U(r-2))$ and the gap is $4r-9$. If $r=4$, then $X=\SU(4,4)/S(U(4)\times U(4))$ and $Y=\SU(3,3)/S(U(3)\times U(3))$. The codimension takes its second minimal value $16$ when $Y'=\SL(4,\mathbb C)/\SU(4)$, hence the gap is $3$. If $r=5$, then $X=\SU(5,5)/S(U(5)\times U(5))$ and $Y=\SU(4,4)/S(U(4)\times U(4))$. The codimension takes its second minimal value $25$ when $Y'$ is $\SL(5,\mathbb C)/\SU(5)$, hence the gap is $8$. If $r=6$, then $X=\SU(6,6)/S(U(6)\times U(6))$ and $Y=\SU(5,5)/S(U(5)\times U(5))$. The codimension takes its second minimal value $36$ when $Y'$ is $\mathbb H^3\times \SU(4,4)/S(U(4)\times U(4))$ or $\SL(6,\mathbb C)/\SU(6)$, hence the gap is $15$.

\emph{Case of $\Sp(r,\mathbb{C})/\Sp(r)$}: the Dynkin diagram is of type $C_r$ and is shown in Figure 3, with multiplicities $2$ for all simple roots. If we remove $\alpha_i$, the remaining diagram will represent $Y_i=\SL(i,\mathbb C)/\SU(i)\times \Sp(r-i,\mathbb{C})/\Sp(r-i))$ (notice $\SL(1,\mathbb C)/\SU(1)$ and $\Sp(0,\mathbb{C})/\Sp(0)$ are just a point by abuse of notation). We compute that the codimension of $Y_i\times \mathbb R$ in $X$ is $-3i^2+(4r+1)i$. It has minimal value $4r-2$ when $i=1$, provided $r>2$. If $r=2$, then $X$ is isomorphic to $\SO(5,\mathbb C)/\SO(5)$, which has been solved previously. Now the codimension takes the second minimal value $8r-10$ when $i=2$, provided $r>5$. In this case, $Y'=\mathbb H^3\times \Sp(r-2,\mathbb{C})/\Sp(r-2)$ and the gap is $4r-8$. If $r=4$, then $X=\Sp(4,\mathbb{C})/\Sp(4)$ and $Y=\Sp(3,\mathbb{C})/\Sp(3)$. The codimension takes its second minimal value $20$ when $Y'=\SL(4,\mathbb C)/\SU(4)$, hence the gap is $6$. If $r=5$, then $X=\Sp(5,\mathbb{C})/\Sp(5)$ and $Y=\Sp(4,\mathbb{C})/\Sp(4)$. The codimension takes its second minimal value $30$ when $Y'$ is $\mathbb H^3\times \Sp(3,\mathbb C)/\Sp(3)$ or $\SL(5,\mathbb C)/\SU(5)$, hence the gap is $12$.

\emph{Case of $\SO^*(4r)/U(2r)$}: the Dynkin diagram is of type $C_r$ and is shown in Figure 3, with ordered multiplicities $4,4,...,4,1$. If we remove $\alpha_i$, the remaining diagram will represent $Y_i=\SU^*(2i-2)/\Sp(i)\times \SO^*(4r-4i)/U(2r-2i)$ (notice $\SU^*(0)/\Sp(1)$ and $\SO^*(0)/U(0)$ are just a point by abuse of notation). We compute that the codimension of $Y_i\times \mathbb R$ in $X$ is $-6i^2+(8r-1)i$. It has minimal value $8r-7$ when $i=1$, provided $r>3$. If $r=2$, then $X$ is isomorphic to $\SO_0(2,6)/\SO(2)\times\SO(6)$, which has been solved previously. If $r=3$, then $X=\SO^*(12)/U(6)$, and the codimension has minimal value $15$ when $i=3$. In this case, the splitting rank occurs when $Y=\SU^*(6)/\Sp(3)$ and is equal to $15$. Now the codimension takes the second minimal value $16r-26$ when $i=2$, provided $r>6$. In this case, $Y'=\mathbb H^5\times \SO^*(4r-8)/U(2r-4)$ and the gap is $8r-19$. If $r=4$, then $X=\SO^*(16)/U(8)$ and $Y=\SO^*(12)/U(6)$. The codimension takes its second minimal value $28$ when $Y'=\SU^*(8)/\Sp(4)$, hence the gap is $3$. If $r=5$, then $X=\SO^*(20)/U(10)$ and $Y=\SO^*(16)/U(8)$. The codimension takes its second minimal value $45$ when $Y'=\SU^*(10)/\Sp(5)$, hence the gap is $12$. If $r=6$, then $X=\SO^*(24)/U(12)$ and $Y=\SO^*(20)/U(10)$. The codimension takes its second minimal value $66$ when $Y'=\SU^*(12)/\Sp(6)$, hence the gap is $25$.

\emph{Case of $\Sp(r,r)/\Sp(r)\times \Sp(r)$}: the Dynkin diagram is of type $C_r$ and is shown in Figure 3, with ordered multiplicities $4,4,...,4,3$. If we remove $\alpha_i$, the remaining diagram will represent $Y_i=\SU^*(2i-2)/\Sp(i)\times [\Sp(r-i,r-i)/\Sp(r-i)\times \Sp(r-i)]$ (Notice $\SU^*(0)/\Sp(1)$ and $\Sp(0,0)/\Sp(0)\times \Sp(0)$ are just a point by abuse of notation). We compute that the codimension of $Y_i\times \mathbb R$ in $X$ is $-6i^2+(8r+1)i$. It has minimal value $8r-5$ when $i=1$, provided $r>2$. If $r=2$, then $X=\Sp(2,2)/\Sp(2)\times \Sp(2)$, and the codimension has minimal value $10$ when $i=2$. So the splitting rank occurs when $Y=\mathbb H^5$ and is equal to $6$. Now the codimension takes the second minimal value $16r-22$ when $i=2$, provided $r>5$. In this case, $Y'=\mathbb H^5\times [\Sp(r-2,r-2)/\Sp(r-2)\times \Sp(r-2)]$ and the gap is $8r-17$. If $r=4$, then $X=\Sp(4,4)/\Sp(4)\times \Sp(4)$ and $Y=\Sp(3,3)/\Sp(3)\times \Sp(3))$. The codimension takes its second minimal value $36$ when $Y'=\SU^*(8)/\Sp(4)$, hence the gap is $9$. If $r=5$, then $X=\Sp(5,5)/\Sp(5)\times \Sp(5)$ and $Y=\Sp(4,4)/\Sp(4)\times \Sp(4)$. The codimension takes its second minimal value $55$ when $Y'=\SU^*(10)/\Sp(5)$, hence the gap is $20$.

\emph{Case of $E_7^{-25}/E_6\times U(1)$}:
the Dynkin diagram is of type $C_3$ and is shown in Figure 3 where $r=3$, with ordered multiplicities $8,8,1$. Hence $Y$ can only be $\SO_0(2,10)/\SO(2)\times \SO(10)$ (when removing $\alpha_1$), or $\mathbb H^9\times \mathbb H^2$ (when removing $\alpha_2$), or $E_6^{-26}/F_4$ (when removing $\alpha_3$). Among the three spaces, $E_6^{-26}/F_4$ has largest dimension thus $\text{srk}(X)=27$. Notice that $E_7^{-25}/E_6\times U(1)$ is of rank three and it does not satisfy the condition of Proposition \ref{prop:gap}.

\begin{figure}
\makebox[1pt][c]{
\beginpicture
\setcoordinatesystem units <.9cm,.9cm> point at -.4 1
\setplotarea x from -.4 to 5, y from -1.4 to 1.4
\linethickness=0.7pt
\putrule from 0.1 0  to 0.9 0
\putrule from 2.1 0  to 2.9 0
\put {$\cdot$} at 1.2 0
\put {$\cdot$} at 1.3 0
\put {$\cdot$} at 1.4 0
\put {$\cdot$} at 1.5 0
\put {$\cdot$} at 1.6 0
\put {$\cdot$} at 1.7 0
\put {$\cdot$} at 1.8 0
\put {$\circ$} at  0  0
\put {$\circ$} at  1 0
\put {$\circ$} at  2  0
\put {$\circ$} at 3  0
\put {$\circ$} at 3.4  0.4
\put {$\circ$} at 3.4  -0.4
\put {$\alpha_1$} at  0 -.3
\put {$\alpha_2$} at  1 -.3
\put {$\alpha_{r-3}$} at  1.8 -.3
\put {$\alpha_{r-2}$} at  2.8 -.3
\put {$\alpha_{r-1}$} at  4 .4
\put {$\alpha_r$} at 4 -.4
\put {$\diagup$} at 3.2 0.2
\put {$\diagdown$} at 3.2 -0.2

\endpicture
}
\caption{Dynkin diagram of type $D_r$}
\end{figure}

\emph{Case of $\SO_0(r,r)/\SO(r)\times \SO(r)$}: the Dynkin diagram is of type $D_r$ and is shown in Figure 4, with multiplicities $1$ for all simple roots. If we remove $\alpha_i$, the remaining diagram will represent $Y_i=\SL(i,\mathbb R)/\SO(i)\times [\SO_0(r-i,r-i)/\SO(r-i)\times \SO(r-i)]$ when $i< r-2$ (notice $\SO_0(3,3)/\SO(3)\times \SO(3)$ is the same as $\SL(4,\mathbb R)/\SO(4)$), and $Y_i=\mathbb H^2\times \mathbb H^2\times\SL(r-2,\mathbb R)/\SO(r-2)$ when $i=r-2$, and $Y_i=\SL(r,\mathbb R)/\SO(r)$ when $i=r-1, r$. We compute that the codimension of $Y_i\times \mathbb R$ in $X$ is $-3i^2/2+(4r-1)i/2$ for $1\leq i\leq r-2$ or $i=r$. It has unique minimal value $2r-2$ when $i=1$, provided $r>4$. If $r=4$, the codimension has minimal value $6$ when $i=1,3,4$. So the splitting rank occurs when $Y=\SL(4,\mathbb R)/\SO(4)$ and is equal to $10$. This agrees with the general result when $r>4$ hence can be absorbed into it. Now the codimension takes the second minimal value $4r-7$ when $i=2$, provided $r>7$. In this case, $Y'=\mathbb H^2\times [\SO_0(r-2,r-2)/\SO(r-2)\times \SO(r-2)]$ and the gap is $2r-5$. If $r=4$, then $X=\SO_0(4,4)/\SO(4)\times \SO(4)$ and $Y=\SL(4,\mathbb R)/\SO(4)$. The codimension takes its second minimal value $9$ when $Y'=\mathbb H^2\times \mathbb H^2\times \mathbb H^2$, hence the gap is $3$. If $r=5$, then $X=\SO_0(5,5)/\SO(5)\times \SO(5)$ and $Y=\SO_0(4,4)/\SO(4)\times \SO(4)$. The codimension takes its second minimal value $10$ when $Y'=\SL(5,\mathbb R)/\SO(5)$, hence the gap is $2$. If $r=6$, then $X=\SO_0(6,6)/\SO(6)\times \SO(6)$ and $Y=\SO_0(5,5)/\SO(5)\times \SO(5)$. The codimension takes its second minimal value $15$ when $Y'=\SL(6,\mathbb R)/\SO(6)$, hence the gap is $5$. If $r=7$, then $X=\SO_0(7,7)/\SO(7)\times \SO(7)$ and $Y=\SO_0(6,6)/\SO(6)\times \SO(6)$. The codimension takes its second minimal value $21$ when $Y'$ is either $\mathbb H^2\times\SO_0(5,5)/\SO(5)\times\SO(5)$ or $\SL(7,\mathbb R)/\SO(7)$, hence the gap is $9$.

\emph{Case of $\SO(2r,\mathbb{C})/\SO(2r)$}: the Dynkin diagram is of type $D_r$ and is shown in Figure 4, with multiplicities $2$ for all simple roots. If we remove $\alpha_i$, the remaining diagram will represent $Y_i=\SL(i,\mathbb C)/\SU(i)\times \SO(2r-2i,\mathbb{C})/\SO(2r-2i)$ when $i< r-2$ (notice $\SO(6,\mathbb{C})/\SO(6)$ is the same as $\SL(4,\mathbb C)/\SU(4)$), and $Y_i=\mathbb H^3\times \mathbb H^3\times\SL(r-2,\mathbb C)/\SU(r-2)$ when $i=r-2$, and $Y_i=\SL(r,\mathbb C)/\SU(r)$ when $i=r-1, r$. We compute that the codimension of $Y_i\times \mathbb R$ in $X$ is $-3i^2+(4r-1)i$ for $1\leq i\leq r-2$ or $i=r$. It has unique minimal value $4r-4$ when $i=1$, provided $r>4$. If $r=4$, the codimension has minimal value $12$ when $i=1,3,4$. So the splitting rank occurs when $Y=\SL(4,\mathbb C)/\SU(4)$ and is equal to $16$. This agrees with the general result when $r>4$ hence can be absorbed into it. Now the codimension takes the second minimal value $8r-14$ when $i=2$, provided $r>7$. In this case, $Y'=\mathbb H^3\times \SO(2r-4,\mathbb{C})/\SO(2r-4)$ and the gap is $4r-10$. If $r=4$, then $X=\SO(8,\mathbb{C})/\SO(8)$ and $Y=\SL(4,\mathbb C)/\SU(4)$. The codimension takes its second minimal value $18$ when $Y'=\mathbb H^3\times \mathbb H^3\times \mathbb H^3$, hence the gap is $6$. If $r=5$, then $X=\SO(10,\mathbb{C})/\SO(10)$ and $Y=\SO(8,\mathbb{C})/\SO(8)$. The codimension takes its second minimal value $20$ when $Y'=\SL(5,\mathbb C)/\SU(5)$, hence the gap is $4$. If $r=6$, then $X=\SO(12,\mathbb{C})/\SO(12)$ and $Y=\SO(10,\mathbb{C})/\SO(10)$. The codimension takes its second minimal value $30$ when $Y'=\SL(6,\mathbb C)/\SU(6)$, hence the gap is $10$. If $r=7$, then $X=\SO(14,\mathbb{C})/\SO(14)$ and $Y=\SO(12,\mathbb{C})/\SO(12)$. The codimension takes its second minimal value $42$ when $Y'$ is either $\mathbb H^3\times \SO(10,\mathbb{C})/\SO(10)$ or $\SL(7,\mathbb C)/\SU(7)$, hence the gap is $18$.

\begin{figure}
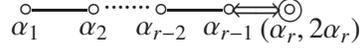

\makebox[1pt][c]{
\beginpicture
\setcoordinatesystem units <.9cm,.9cm> point at -.4 1
\setplotarea x from -.4 to 4.4, y from -1.4 to 1.4
\linethickness=.7pt
\putrule from 0.1 0  to 0.9 0
\putrule from 2.1 0  to 2.9 0
\put {$\cdot$} at 1.2 0
\put {$\cdot$} at 1.3 0
\put {$\cdot$} at 1.4 0
\put {$\cdot$} at 1.5 0
\put {$\cdot$} at 1.6 0
\put {$\cdot$} at 1.7 0
\put {$\cdot$} at 1.8 0
\put {$\circ$} at  0  0
\put {$\circ$} at  1 0
\put {$\circ$} at  2  0
\put {$\circ$} at 3  0
\put {$\circ$} at 3.9  0
\put {$\bigcirc$} at 3.9 0
\put {$\Longleftrightarrow$} at 3.4 0
\put {$\alpha_1$} at  0 -.3
\put {$\alpha_2$} at  1 -.3
\put {$\alpha_{r-2}$} at  2 -.3
\put {$\alpha_{r-1}$} at  3 -.3
\put {$(\alpha_r,2\alpha_r)$} at 4.2 -.3
\endpicture
}
\caption{Dynkin diagram of type $(BC)_r$}
\end{figure}

\emph{Case of $ \SU(r,r+k)/S(U(r)\times U(r+k))$}: the Dynkin diagram is of type $(BC)_r$ and is shown in Figure 5, with ordered multiplicities $2,2,...,2,(2k,1)$. If we remove $\alpha_i$, the remaining diagram will represent $Y_i=\SL(i,\mathbb C)/\SU(i)\times \SU(r-i,r-i+k)/S(U(r-i)\times U(r-i+k))$ (notice $\SL(1,\mathbb C)/\SU(1)$ and $\SU(0,k)/S(U(0)\times U(k))$ are just a point by abuse of notation). We compute that the codimension of $Y_i\times \mathbb R$ in $X$ is $-3i^2+(4r+2k)i$. It has unique minimal value $4r+2k-3$ when $i=1$, provided $r+2k>3$, which holds for higher rank symmetric spaces. So the splitting rank occurs when $Y=\SU(r-1,r-1+k)/S(U(r-1)\times U(r-1+k))$ and is equal to $n-(4r+2k-3)$. Now the codimension takes the second minimal value $8r+4k-12$ when $i=2$, provided $r+2k>6$. In this case, $Y'=\mathbb H^3\times\SU(r-2,r-2+k)/S(U(r-2)\times U(r-2+k))$ and the gap is $4r+2k-9$. As we focus on $r\geq 4$ in Proposition \ref{prop:gap}, the only space excluded by $r+2k>6$ is $\SU(4,5)/S(U(4)\times U(5))$ ($r=4,k=1$). In this case, $Y$ is $\SU(3,4)/S(U(3)\times U(4))$. The codimension takes its second minimal value $24$ when $Y'$ is either $\mathbb H^3\times \SU(2,3)/S(U(2)\times U(3))$ or $\SL(4,\mathbb C)/\SU(4)$, hence the gap is $9$.

\emph{Case of $\Sp(r,r+k)/\Sp(r)\times \Sp(r+k)$}: the Dynkin diagram is of type $(BC)_r$ and is shown in Figure 5, with ordered multiplicities $4,4,...,4,(4k,3)$. If we remove $\alpha_i$, the remaining diagram will represent $Y_i=\SU^*(2i)/\Sp(i)\times [\Sp(r-i,r-i+k)/\Sp(r-i)\times \Sp(r-i+k)]$ (notice $\SU^*(2)/\Sp(1)$ and $\Sp(0,k)/\Sp(0)\times \Sp(k)$ are just a point by abuse of notation). We compute that the codimension of $Y_i\times \mathbb R$ in $X$ is $-6i^2+(8r+4k+1)i$. It has unique minimal value $8r+4k-5$ when $i=1$. So the splitting rank occurs when $Y=\Sp(r-1,r-1+k)/\Sp(r-1)\times \Sp(r-1+k)$ and is equal to $n-(8r+4k-5)$. Now the codimension takes the second minimal value $16r+8k-22$ when $i=2$, provided $2r+4k>11$. In this case, $Y'=\mathbb H^5\times\Sp(r-2,r-2+k)/\Sp(r-2)\times \Sp(r-2+k)$ and the gap is $8r+4k-17$. As we focus on $r\geq 4$ in Proposition \ref{prop:gap}, the inequality $2r+4k>11$ always holds.

\emph{Case of $\SO^*(4r+2)/U(2r+1)$}: the Dynkin diagram is of type $(BC)_r$ and is shown in Figure 5, with ordered multiplicities $4,4,...,4,(4,1)$. If we remove $\alpha_i$, the remaining diagram will represent $Y_i=\SU^*(2i)/\Sp(i)\times \SO^*(4r-4i+2)/U(2r-2i+1)$ (notice $\SU^*(2)/\Sp(1)$ and $\SO^*(2)/U(1)$ are just a point by abuse of notation).  We compute that the codimension of $Y_i\times \mathbb R$ in $X$ is $-6i^2+(8r+3)i$. It has unique minimal value $8r-3$ when $i=1$. So the splitting rank occurs when $Y=\SO^*(4r-2)/U(2r-1)$ and is equal to $n-(8r-3)$.  Now the codimension takes the second minimal value $16r-18$ when $i=2$, provided $r>4$. In this case, $Y'=\mathbb H^5\times\SO^*(4r-6)/U(2r-3)$ and the gap is $8r-15$. As we focus on $r\geq 4$ in Proposition \ref{prop:gap}, the only space excluded by the inequality $r>4$ is $\SO^*(18)/U(9)$. In this special case, $Y=\SO^*(14)/U(7)$ and the codimension takes its second minimal value $44$ when $Y'=\SU^*(8)/\Sp(4)$, hence the gap is $15$.

\emph{Case of $E_6^{-14}/\text{Spin}(10)\times U(1)$}: the Dynkin diagram is of type $(BC)_2$ and is shown in Figure 5, with ordered multiplicities $6,(8,1)$. Hence $Y$ can only be $\mathbb H^9$ or $\SU(1,5)/S(U(1)\times U(5))\simeq \mathbb C\mathbb H^5$. Comparing the dimensions of the two spaces, we conclude the one that has splitting rank should be $\mathbb C\mathbb H^5\times \mathbb R$, and the splitting rank is $11$. Notice that $E_6^{-14}/\text{Spin}(10)\times U(1)$ is of rank two hence it does not satisfy the condition of Proposition \ref{prop:gap}.

\begin{figure}
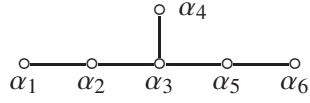

\makebox[1pt][c]{
\beginpicture
\setcoordinatesystem units <.9cm,.9cm> point at -.4 1
\setplotarea x from -.4 to 4.4, y from -1.4 to 1.4
\linethickness=.7pt
\putrule from 0.1 0  to 0.9 0
\putrule from 2.1 0  to 2.9 0
\putrule from 1.1 0  to 1.9 0
\putrule from 3.1 0  to 3.9 0
\putrule from 2 0.1  to 2 0.7
\put {$\circ$} at  0  0
\put {$\circ$} at  1 0
\put {$\circ$} at  2  0
\put {$\circ$} at 3  0
\put {$\circ$} at 4  0
\put {$\circ$} at 2  0.8
\put {$\alpha_1$} at  0 -.3
\put {$\alpha_2$} at  1 -.3
\put {$\alpha_3$} at  2 -.3
\put {$\alpha_4$} at  2.5 0.8
\put {$\alpha_5$} at  3 -.3
\put {$\alpha_6$} at  4 -.3
\endpicture
}
\caption{Dynkin diagram of type $E_6$}
\end{figure}

\emph{Case of $E_6^6/\Sp(4)$}: the Dynkin diagram is of type $E_6$ and is shown in Figure 6, with multiplicities $1$ for all simple roots. If we remove one simple root, the remaining diagram will represent $4$ kinds of symmetric spaces: $Y_1=Y_6=\SO_0(5,5)/\SO(5)\times \SO(5)$, $Y_2=Y_5=\mathbb H^2\times \SL(5,\mathbb R)/\SO(5)$, $Y_3=\mathbb H^2\times \SL(3,\mathbb R)/SO(3)\times \SL(3,\mathbb R)/\SO(3)$ and $Y_4=\SL(6,\mathbb R)/\SO(6)$. We compute that the dimensions of $Y_i\times \mathbb R$ are $26,17,13$ and $21$ respectively. So the splitting rank is $26$ and the gap is $5$.

\emph{Case of $E_6(\mathbb C)/E_6$}: the Dynkin diagram is of type $E_6$ and is shown in Figure 6, with multiplicities $2$ for all simple roots. If we remove one simple root, the remaining diagram will represent $4$ kinds of symmetric spaces: $Y_1=Y_6=\SO(10,\mathbb C)/\SO(10)$, $Y_2=Y_5=\mathbb H^3\times \SL(5,\mathbb C)/\SU(5)$, $Y_3=\mathbb H^3\times \SL(3,\mathbb C)/\SU(3)\times \SL(3,\mathbb C)/\SU(3)$ and $Y_4=\SL(6,\mathbb C)/\SU(6)$. We compute that the dimensions of $Y_i\times \mathbb R$ are $46,28,20$ and $36$ respectively. So the splitting rank is $46$ and the gap is $10$.

\begin{figure}
\makebox[1pt][c]{
\beginpicture
\setcoordinatesystem units <.9cm,.9cm> point at -.4 1
\setplotarea x from -.4 to 5.4, y from -1.4 to 1.4
\linethickness=.7pt
\putrule from 0.1 0  to 0.9 0
\putrule from 2.1 0  to 2.9 0
\putrule from 1.1 0  to 1.9 0
\putrule from 3.1 0  to 3.9 0
\putrule from 2 0.1  to 2 0.7
\putrule from 4.1 0  to 4.9 0
\put {$\circ$} at  0  0
\put {$\circ$} at  1 0
\put {$\circ$} at  2  0
\put {$\circ$} at 3  0
\put {$\circ$} at 4  0
\put {$\circ$} at 5  0
\put {$\circ$} at 2  0.8
\put {$\alpha_1$} at  0 -.3
\put {$\alpha_2$} at  1 -.3
\put {$\alpha_3$} at  2 -.3
\put {$\alpha_4$} at  2.5 0.8
\put {$\alpha_5$} at  3 -.3
\put {$\alpha_6$} at  4 -.3
\put {$\alpha_7$} at  5 -.3
\endpicture
}
\caption{Dynkin diagram of type $E_7$}
\end{figure}

\emph{Case of $E_7^7/\SU(8)$}: the Dynkin diagram is of type $E_7$ and is shown in Figure 7, with multiplicities $1$ for all simple roots. If we remove one simple root, the remaining diagram will represent $7$ kinds of symmetric spaces: $Y_1=\SO_0(6,6)/\SO(6)\times \SO(6)$, $Y_2=\mathbb H^2\times \SL(6,\mathbb R)/\SO(6)$, $Y_3=\mathbb H^2\times \SL(3,\mathbb R)/\SO(3)\times \SL(4,\mathbb R)/\SO(4)$, $Y_4=\SL(7,\mathbb R)/\SO(7)$, $Y_5=\SL(3,\mathbb R)/\SO(3)\times \SL(5,\mathbb R)/\SO(5)$, $Y_6=\mathbb H^2\times \SO_0(5,5)/\SO(5)\times \SO(5)$ and $Y_7=E_6^6/\Sp(4)$. We compute that the dimensions of $Y_i\times \mathbb R$ are $37,23,17,28,20,28$ and $43$ respectively. So the splitting rank is $43$ and the gap is $6$.

\emph{Case of $E_7(\mathbb C)/E_7$}: the Dynkin diagram is of type $E_7$ and is shown in Figure 7, with multiplicities $2$ for all simple roots. If we remove one simple root, the remaining diagram will represent $7$ kinds of symmetric spaces: $Y_1=\SO(12,\mathbb C)/\SO(12)$, $Y_2=\mathbb H^3\times \SL(6,\mathbb C)/\SU(6)$, $Y_3=\mathbb H^3\times \SL(3,\mathbb C)/\SU(3)\times \SL(4,\mathbb C)/\SU(4)$, $Y_4=\SL(7,\mathbb C)/\SU(7)$, $Y_5=\SL(3,\mathbb C)/\SU(3)\times \SL(5,\mathbb C)/\SU(5)$, $Y_6=\mathbb H^3\times \SO(10,\mathbb C)/\SO(10)$ and $Y_7=E_6(\mathbb C)/E_6$. We compute that the dimensions of $Y_i\times \mathbb R$ are $67,39,27,49,33,49$ and $79$ respectively. So the splitting rank is $79$ and the gap is $12$.

\begin{figure}
\makebox[1pt][c]{
\beginpicture
\setcoordinatesystem units <.9cm,.9cm> point at -.4 1
\setplotarea x from -.4 to 6.4, y from -1.4 to 1.4
\linethickness=.7pt
\putrule from 0.1 0  to 0.9 0
\putrule from 2.1 0  to 2.9 0
\putrule from 1.1 0  to 1.9 0
\putrule from 3.1 0  to 3.9 0
\putrule from 2 0.1  to 2 0.7
\putrule from 4.1 0  to 4.9 0
\putrule from 5.1 0  to 5.9 0

\put {$\circ$} at  0  0
\put {$\circ$} at  1 0
\put {$\circ$} at  2  0
\put {$\circ$} at 3  0
\put {$\circ$} at 4  0
\put {$\circ$} at 5  0
\put {$\circ$} at 6  0
\put {$\circ$} at 2  0.8
\put {$\alpha_1$} at  0 -.3
\put {$\alpha_2$} at  1 -.3
\put {$\alpha_3$} at  2 -.3
\put {$\alpha_4$} at  2.5 0.8
\put {$\alpha_5$} at  3 -.3
\put {$\alpha_6$} at  4 -.3
\put {$\alpha_7$} at  5 -.3
\put {$\alpha_8$} at  6 -.3
\endpicture
}
\caption{Dynkin diagram of type $E_8$}
\end{figure}

\emph{Case of $E_8^8/\SO(16)$}: the Dynkin diagram is of type $E_8$ and is shown in Figure 8, with multiplicities $1$ for all simple roots. If we remove one simple root, the remaining diagram will represent $8$ kinds of symmetric spaces: $Y_1=\SO_0(7,7)/\SO(7)\times \SO(7)$, $Y_2=\mathbb H^2\times \SL(7,\mathbb R)/\SO(7)$, $Y_3=\mathbb H^2\times \SL(3,\mathbb R)/\SO(3)\times \SL(5,\mathbb R)/\SO(5)$, $Y_4=\SL(8,\mathbb R)/\SO(8)$, $Y_5=\SL(4,\mathbb R)/\SO(4)\times \SL(5,\mathbb R)/\SO(5)$, $Y_6=\SL(3,\mathbb R)/\SO(3)\times \SO_0(5,5)/\SO(5)\times \SO(5)$, $Y_7=\mathbb H^2\times E_6^6/\Sp(4)$, and $Y_8=E_7^7/\SU(8)$. We compute that the dimensions of $Y_i\times \mathbb R$ are $50,30,22,36,24,31,45$ and $71$ respectively. So the splitting rank is $71$ and the gap is $21$.

\emph{Case of $E_8(\mathbb C)/E_8$}: the Dynkin diagram is of type $E_8$ and is shown in Figure 8, with multiplicities $2$ for all simple roots. If we remove one simple root, the remaining diagram will represent $8$ kinds of symmetric spaces: $Y_1=\SO(14,\mathbb C)/\SO(14)$, $Y_2=\mathbb H^3\times \SL(7,\mathbb C)/\SU(7)$, $Y_3=\mathbb H^3\times \SL(3,\mathbb C)/\SU(3)\times \SL(5,\mathbb C)/\SU(5)$, $Y_4=\SL(8,\mathbb C)/\SU(8)$, $Y_5=\SL(4,\mathbb C)/\SU(4)\times \SL(5,\mathbb C)/\SU(5)$, $Y_6=\SL(3,\mathbb C)/\SU(3)\times \SO(10,\mathbb C)/\SO(10)$, $Y_7=\mathbb H^3\times E_6(\mathbb C)/E_6)$, and $Y_8=E_7(\mathbb C)/E_7$. We compute that the dimensions of $Y_i\times \mathbb R$ are $92,52,36,64,40,54,82$ and $134$ respectively. So the splitting rank is $134$ and the gap is $42$.

\begin{figure}
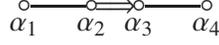

\makebox[1pt][c]{
\beginpicture
\setcoordinatesystem units <.9cm,.9cm> point at -.4 1
\setplotarea x from -.4 to 3.4, y from -1.4 to 1.4
\linethickness=.7pt
\putrule from 0.1 0  to 0.9 0
\putrule from 1.8 0  to 2.6 0
\put {$\Longrightarrow$} at 1.35 0
\put {$\circ$} at  0  0
\put {$\circ$} at  1 0
\put {$\circ$} at 1.7  0
\put {$\circ$} at 2.7 0
\put {$\alpha_1$} at  0 -.3
\put {$\alpha_2$} at  1 -.3
\put {$\alpha_3$} at  1.7 -.3
\put {$\alpha_4$} at  2.7 -.3
\endpicture
}
\caption{Dynkin diagram of type $F_4$}
\end{figure}

\emph{Case of $F_4^4/\Sp(3)\times\Sp(1)$}: the Dynkin diagram is of type $F_4$ and is shown in Figure 9, with multiplicities $1$ for all simple roots. If we remove one simple root, the remaining diagram will represent $3$ kinds of symmetric spaces: $Y_1=\Sp(3,\mathbb R)/U(3)$, $Y_2=Y_3=\mathbb H^2\times \SL(3,\mathbb R)/\SO(3)$, and $Y_4=\SO_0(3,4)/\SO(3)\times \SO(4)$. We compute that the dimensions of $Y_i\times \mathbb R$ are $13, 8$ and $13$ respectively. So the splitting rank is $13$ and the gap is $5$.

\emph{Case of $E_6^2/\SU(6)\times\Sp(1)$}: the Dynkin diagram is of type $F_4$ and is shown in Figure 9, with ordered multiplicities $1,1,2,2$. If we remove one simple root, the remaining diagram will represent $4$ kinds of symmetric spaces: $Y_1=\SU(3,3)/S(U(3)\times U(3))$, $Y_2=\mathbb H^2\times \SL(3,\mathbb C)/\SU(3)$, $Y_3=\mathbb H^3\times \SL(3,\mathbb R)/\SO(3)$, and $Y_4=\SO_0(3,5)/\SO(3)\times \SO(5)$. We compute that the dimensions of $Y_i\times \mathbb R$ are $19, 11, 9$ and $16$ respectively. So the splitting rank is $19$ and the gap is $3$.

\emph{Case of $E_7^{-5}/\SO(12)\times\Sp(1)$}: the Dynkin diagram is of type $F_4$ and is shown in Figure 9, with ordered multiplicities $1,1,4,4$. If we remove one simple root, the remaining diagram will represent $4$ kinds of symmetric spaces: $Y_1=\SO^*(12)/U(6)$, $Y_2=\mathbb H^2\times \SU^*(6)/\Sp(3)$, $Y_3=\mathbb H^5\times \SL(3,\mathbb R)/\SO(3)$, and $Y_4=\SO_0(3,7)/\SO(3)\times \SO(7)$. We compute that the dimensions of $Y_i\times \mathbb R$ are $31, 17, 11$ and $28$ respectively. So the splitting rank is $31$ and the gap is $3$.

\emph{Case of $E_8^{-24}/E_7\times \Sp(1)$}: the Dynkin diagram is of type $F_4$ and is shown in Figure 9, with ordered multiplicities $1,1,8,8$. If we remove one simple root, the remaining diagram will represent $4$ kinds of symmetric spaces: $Y_1=E_7^{-25}/E_6\times U(1)$, $Y_2=\mathbb H^2\times E_6^{-26}/F_4$, $Y_3=\mathbb H^9\times \SL(3,\mathbb R)/\SO(3)$, and $Y_4=\SO_0(3,11)/\SO(3)\times \SO(11)$. We compute that the dimensions of $Y_i\times \mathbb R$ are $55, 29, 15$ and $34$ respectively. So the splitting rank is $55$ and the gap is $21$.

\emph{Case of $F_4(\mathbb C)/F_4$}: the Dynkin diagram is of type $F_4$ and is shown in Figure 9, with multiplicities $2$ for all simple roots. If we remove one simple root, the remaining diagram will represent $3$ kinds of symmetric spaces: $Y_1=\Sp(3,\mathbb C)/\Sp(3)$, $Y_2=Y_3=\mathbb H^3\times \SL(3,\mathbb C)/\SU(3)$, and $Y_4=\SO(7,\mathbb C)/\SO(7)$. We compute that the dimensions of $Y_i\times \mathbb R$ are $22, 12$ and $22$ respectively. So the splitting rank is $22$ and the gap is $10$.

\begin{figure}
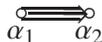

\makebox[1pt][c]{
\beginpicture
\setcoordinatesystem units <.9cm,.9cm> point at -.4 1
\setplotarea x from -.4 to 1.4, y from -1.4 to 1.4
\linethickness=.7pt
\putrule from 0 0.09  to 0.85 0.09
\putrule from 0 -0.07  to 0.85 -0.07
\putrule from 0.1 0.01  to 0.9 0.01
\put {$\circ$} at  0  0
\put {$\circ$} at  1 0
\put {$\alpha_1$} at  0 -.3
\put {$\alpha_2$} at  1 -.3
\put {$\Rrightarrow$} at 0.8 0
\endpicture
}
\caption{Dynkin diagram of type $G_2$}
\end{figure}

\emph{Case of $G_2^2/\SO(4)$}: the Dynkin diagram is of type $G_2$ and is shown in Figure 10, with multiplicities $1$ for both simple roots. If we remove one simple root, the remaining diagram will represent the only symmetric space: $\mathbb H^2$. So the splitting rank is $3$ corresponding to the totally geodesic submanifold $\mathbb H^2\times \mathbb R$. Notice this space is of rank two so it does not satisfy the condition of Proposition \ref{prop:gap}.

\emph{Case of $G_2(\mathbb C)/G_2$}: the Dynkin diagram is of type $G_2$ and is shown in Figure 10, with multiplicities $2$ for both simple roots. If we remove one simple root, the remaining diagram will represent the only symmetric space: $\mathbb H^3$. So the splitting rank is $4$ corresponding to the totally geodesic submanifold $\mathbb H^3\times \mathbb R$. Notice this space is of rank two so it does not satisfy the condition of Proposition \ref{prop:gap}.

This verifies all cases, and completes the proofs of both Theorem \ref{thm:srk} and Proposition \ref{prop:gap}.

\bibliographystyle{plain}

\end{document}